\documentclass{amsart}
\usepackage{amssymb,latexsym}
\usepackage{graphicx}
\newtheorem{theorem}{Theorem}
\newtheorem{lemma}{Lemma}
\newtheorem{definition}{Definition}

\newtheorem{corollary}{Corollary}
\newtheorem{remark}{Remark}
\usepackage[latin1]{inputenc}
\usepackage[T1]{fontenc}
\usepackage{setspace}
\usepackage{color}
\usepackage{hyperref}

\title[ It\^o formula for Mild Solutions]{\bf { It\^o formula for Mild Solutions of SPDEs with Gaussian and non-Gaussian noise and applications to stability properties}}

 \author{S. Albeverio}
\address{Sergio Albeverio: Institut f\"ur Angewandte Mathematik Abteilung Wahrscheinlichkeitstheorie, 
 Rheinische Friedrich-Wilhelms-Universit\"at Bonn, Endenicher Allee 60, 53115 Bonn, Germany}
 \email{albeverio@uni-bonn.de / albeverio@yahoo.com}
 
 \author{L. Gawarecki}
 \address{Leszek Gawarecki: Department of Mathematics Kettering University, 1700 University Ave Flint, MI 48504 USA}
 \email{lgawarec@kettering.edu}
 
 \author{V. Mandrekar}
 \address{Vidyadhar Mandrekar: Department of Statistics and Probability, Michigan State University, East Lansing, MI 48823 USA}
 \email{atmah@stt.msu.edu / atma1m@gmail.com}
 
 \author{B. R\"udiger}
 \address{Barbara R\"udiger: Bergische Universit\"at Wuppertal, Fakult\"at f\"ur Mathematik und Naturwissenschaften, Gaussstrasse 20, 
 42119 Wuppertal,
 Germany}
 \email{ruediger@uni-wuppertal.de}
 
 \author{B. Sarkar}
 \address{Barun Sarkar: Bergische Universit\"at Wuppertal, Fakult\"at f\"ur Mathematik und Naturwissenschaften, Gaussstrasse 20, 
 42119 Wuppertal,
 Germany}
 \email{barunsarkar.math@gmail.com}



\begin{document}







 \begin{abstract}
We use  Yosida approximation to find an It\^o  formula for mild solutions $\left\{X^x(t), t\geq 0\right\}$ of SPDEs with  Gaussian and non-Gaussian coloured  noise, the non Gaussian noise being defined through compensated Poisson random 
measure associated to a L\'evy process. The functions to which we apply such  It\^o formula are in $C^{1,2}([0,T]\times H)$, as in the case considered for SDEs in \cite{9}. Using this It\^o formula we prove exponential stability and exponential ultimate boundedness properties in mean square sense for mild solutions.
We also compare such It\^o formula to an  It\^o formula  for mild solutions introduced by Ichikawa in \cite{8}, and an It\^o formula written in terms of the semigroup of the drift operator \cite{11} which we extend before to the non Gaussian case.   
\end{abstract}

\maketitle


\section{Introduction}\label{intro}

 The It\^o formula for  strong  
 solutions of SPDEs can be derived similarly as for the case of  SDEs, see e.g. \cite{1}, \cite{rvu-yor} for the Gaussian case and \cite{9}, \cite{bookR}.  Here we derive through Yosida 
 approximation an It\^o formula for mild solutions $\left\{X^x(t), t\geq 0\right\}$  for SPDEs driven by a Wiener process and general L\'evy 
 processes. The functions to which we apply such  It\^o formula are in $C^{1,2}([0,T]\times H)$, as in the case considered for SDEs in \cite{9}.
 Using this It\^o formula we prove exponential stability and exponential ultimate boundedness properties in mean square sense for mild solutions. We also relate this It\^o formula to an It\^o formula for mild solutions  of stochastic semilinear evolution equations with Gaussian noise provided by Ichikawa in \cite{8}, which we generalize to SPDEs with 
 L\'evy noise. The functions to which such Ichikawa's  It\^o formula can be applied are in the domain of the "weak generator" associated to the SPDE, defined in \eqref{generatorDcor} below.  Following Ichikawa \cite{8} we also show in Section \ref{relt-ichi-sem-gen} that the It\^o formula  for mild solutions $\left\{X^x(t), t\geq 0\right\}$  induces an inequality for $\Psi(X(t))$ which can be applied to a larger set of functions $\Psi$, than those being in the domain of the "weak generator".  Such inequality  might also  be applied to prove exponential stability of $\left\{X^x(t), t\geq 0\right\}$ in mean square sense, as shown through an example.
 
 We also generalize the It\^o formula obtained  in \cite{11} in terms of the semigroup of the drift operator to SPDEs with non Gaussian noise. In \cite{11}  Da Prato, Jentzen and R\"ockner transform the mild solution of a SPDE with Gaussian noise 
  to a standard 
  It\^o process by using the techniques of \cite{14}, through the works of Nagy \cite{nagy1}, 
  \cite{nagy2}, \cite{nagy3}. Different than the  It\^o formula derived in \cite{11} the  It\^o formula obtained here  by Yosida approximation allows us to study the stability and ultimate boundedness properties of the mild solutions of SPDES (see also 
   e.g. \cite{3} for the non-Gaussian case and  \cite{1} for the Gaussian case.)  \\

 In Section \ref{pre}, we give Preliminaries related to the  
  existence and uniqueness of mild solutions and discuss under what assumptions  
the mild solution is also a strong solution.   We recall  how  the mild solutions can be obtained in the limit  by a sequence of strong solutions by using Yosida approximation 
techniques.
\par In 
Section \ref{mthm}  we  recall the It\^o formula for strong solutions. We then derive by Yosida approximation  an  It\^o formula for mild solutions in Theorem \ref{imldtm4}, Section 3.3.  The functions for which we compute the It\^o formula are in  $ C^{1,2}([0,T]\times H)$,  and satisfy the same conditions as those applied  in \cite{9} to strong solutions of SDEs with non Gaussian noise. 
In Section \ref{Applications stability} we prove, using the It\^o  formula obtained in Theorem \ref{imldtm4}, that under suitable conditions the mild solutions of SPDEs are exponential stable (Theorem \ref{m.s.s-th1}) and ultimate bounded (Theorem \ref{ulti-bn-th1}) in mean square sense.

\par 
In \cite{8} Ichikawa obtained an It\^o formula for mild solutions w.r.t. the Gaussian noise using their "weak generators" . This will be recalled in Section 
\ref{ichikawa}. A similar  result will be   derived also for the non-Gaussian noise in Theorem \ref{ichiIto}, Section \ref{ichikawa}. The functions for which we compute the It\^o formula following Ichikawa  are functions $\Psi\in$  $ C^{1,2}([0,T]\times H)$,  
for which the function $\mathcal{L}\Psi$, with  $\mathcal{L}$ being the "weak generator"  \eqref{defnL}, can be extended to a continuous function. Ichikawa's It\^o formula in Theorem \ref{ichiIto} is compared with the It\^o formula obtained in Theorem  \ref{imldtm4} in Section    \ref{ichikawa} through Example 2 and Example 3.
\par
In Ichikawa's type It\^o formula for mild solutions, 
the assumption that 
  the function $\mathcal{L}\Psi$ can be extended to a continuous function, is rather restrictive to study the stability theory of 
 the mild solutions.  
So, following Ichikawa's result for the Gaussian case, in Lemma \ref{newV} of Section \ref{relt-ichi-sem-gen}
we assume that the function $\mathcal{L}\Psi\leq\mathcal{U}$, where 
$\mathcal{U}$ is a continuous function.  
This might be used to study the exponential stability of the mild solutions in the mean square sense (this is shown in Example 4). 
In Section \ref{pratojenroc} we  extend 
Da Prato, Jentzen and R\"ockner's
mild It\^o formula to the case of  non-Gaussian noise in Theorem \ref{TheoremA}. For this  we used the transformation technique presented  in \cite{14}. 
Example 3 in Section 3 can be also obtained through the It\^o formula in Theorem \ref{TheoremA}.

\section{Preliminaries}\label{pre}

Let $K$ and $H$ be real separable Hilbert spaces. Let $(H\backslash\left\{0\right\},\mathcal{B}(H\backslash\left\{0\right\}),\beta)$ be a $\sigma$-finite 
measurable space, with $\mathcal{B}(H\backslash\left\{0\right\})$ denoting the Borel sets of $H\backslash\left\{0\right\}$ and $\beta$ be 
a positive measure on 
$\mathcal{B}(H\backslash\left\{0\right\})$ with

\begin{align}
\int_{H\backslash\left\{0\right\}}(\left\|u\right\|_{H}^{2}\wedge1)\beta(du)<\infty.\nonumber
\end{align}
We refer to this $\beta$ as a L\'evy measure on $H\backslash\left\{0\right\}$.\\

 We shall consider  a compensated Poisson random measure (cPrm)  $q(ds,du):=N(ds,du)(\omega)-ds\beta(du)$ on a filtered probability space 
 $(\Omega,\mathcal{F},\left\{\mathcal{F}_{t}\right\}_{t\in [0, T]},P)$ satisfying the usual hypothesis. 
  $ds$ denotes the Lebesgue measure on $\mathcal{B}(\mathbb{R}_{+})$ and 
 $N(ds,du)(\omega)$ is a Poisson distributed $\sigma$-finite measure on the $\sigma$-algebra $\mathcal{B}(\mathbb{R}_{+},H\backslash\left\{0\right\})$, 
 generated by the product semiring $\mathcal{B}(\mathbb{R}_{+})\times\mathcal{B}(H\backslash\left\{0\right\})$ of the Borel $\sigma$-algebra 
 $\mathcal{B}(\mathbb{R}_{+})$ and the  Borel $\sigma$-algebra $\mathcal{B}(H\backslash\left\{0\right\})$. 
 Then $E(q(A\times B))^{2}=\beta(A)\lambda(B)$, for any $A\in\mathcal{B}(H\backslash\left\{0\right\})$,  $B\in\mathcal{B}(\mathbb{R}_{+})$, $0\notin\bar{A}$, ($\bar{A}$ denoting the closure of the set $A$) 
  $\lambda(B)$ is the Lebesgue measure of $B$. (For more details we refer to section 1 of \cite{4}).\\
 
 The definition of stochastic integral with respect to compensated Poisson random measure and their properties are given in, e.g.  
 \cite{l-ito-dc}, \cite{bookR}, \cite{4}, \cite{6}, \cite{k17}, \cite{2}, \cite{10}, \cite{iki-wat}, \cite{skro}.\\
 We recall here some definition and known facts:
  \begin{definition}\label{xtra-sem-dfn}
   A family $S(t)\in\mathcal{L}(X)$, $t\geq0$, of bounded linear operators on a Banach space $X$ is called a strongly continuous semigroup (or a 
   $C_0$-semigroup) if\\
   (S1) $S(0)=I$,\\
   (S2) (Semigroup property) $S(t+s)=S(t)S(s)$ for every $t,s\geq0$,\\
   (S3) (Strong continuity property) $\lim_{t\rightarrow0^+}S(t)x=x$ for every $x\in X$.
   \end{definition}
   Let $S(t)$ be a $C_0$-semigroup on a Banach space $X$. Then there exist constants $\alpha\geq0$ and $M\geq1$ 
   such that $\|S(t)\|_{\mathcal{L}(X)}\leq Me^{\alpha t}$, $t\geq0$.\\
   
   If $M=1$, then $S(t)$ is called a \textit{pseudo-contraction semigroup}.\\

We consider the following stochastic partial differential equation with values in $H$,

\begin{equation} 
\left\{
\begin{aligned}
dX(t)=(AX(t)+F(X(t)))dt+B(X(t))dW_{t}+\int_{H\backslash\left\{0\right\}} f(v,X(t))q(dv,dt);\\
X(0)=\xi.
\end{aligned}
\right. \label {imldeq1}
\end{equation}

 Where $\xi$ is an $\mathcal{F}_{0}$-measurable random variable. We assume that, the terms in (\ref{imldeq1}) satisfy the following conditions:  \\

(A1) $A$ is the infinitesimal generator of a pseudo-contraction semigroup $\{S(t),t\geq0\}$ on $H$. This means in particular that there exists a constant 
 $\alpha\in\mathbb{R}_+$ s.t. $\left\|S(t)\right\|\leq e^{\alpha t}$. $\mathcal{D}(A)$ denotes the domain of the linear operator $A$. 
 $\mathcal{D}(A)$ is dense in $H$ and $A$ is a closed linear operator. \\

(A2) $(W_t)_{t\geq0}$ is a $K$-valued $\mathcal{F}_{t}$-Wiener process with covariance $Q$ on 
$\left(\Omega,\mathcal{F},\left\{\mathcal{F}_{t}\right\}_{t\in [0,T]},P\right)$ satisfying the usual hypothesis, where $Q$ is a nonnegative definite 
symmetric trace-class operator on the real separable Hilbert space $K$. 
$(W_t)_{t\geq0}$ is assumed to be independent of the cPrm $q(dv,dt)$.    \\

(A3) $F:H\rightarrow H$, $B:H\rightarrow \mathcal{L}(K,H)$, $f:H\backslash\left\{0\right\}\times H\rightarrow H$ are continuous, and Bochner measurable functions satisfying:

\begin{align}
\|F(x)\|_{H}^{2}+tr(B(x)QB^{*}(x))+\int_{H\backslash\left\{0\right\}}\left\|f(v,x)\right\|_{H}^{2}\beta(dv)\leq l(1+\left\|x\right\|_{H}^{2});\nonumber
\end{align}
and
\begin{align*}
\big\|F(x)-F(y)\big\|_{H}^{2}+tr((B(x)-B(y))Q(B(x)-B(y))^{*})
\end{align*}
\begin{align*}
+\int_{H\backslash\left\{0\right\}}\big\|f(v,x)-f(v,y)\big\|_{H}^{2}\beta(dv)\leq \mathcal{K}\big\|x-y\big\|_{H}^{2};
\end{align*}

for all $x,y\in H$. Where $l$, $\mathcal{K}$ are positive constants.  \\

 \begin{definition}\label{imlddf1}
 A stochastic process $\left\{X(t),t\geq0\right\}$ is called a mild solution of (\ref{imldeq1}) in $[0,T]$ \\
 (i) $X(t)$ is $\mathcal{F}_{t}$-adapted on the filtered probability space  ($\Omega,\mathcal{F},\left\{\mathcal{F}_{t}\right\}_{t\in [0,T]},P$), \\
 (ii) $\left\{X(t),t\geq 0\right\}$ is jointly measurable and $\int_{0}^{T}E[\left\|X(t)\right\|_{H}^{2}]dt<\infty$, \\
 (iii)
 \begin{align*}
 X(t)=S(t)\xi+\int_{0}^{t}S(t-s)F(X(s))ds+\int_{0}^{t}S(t-s)B(X(s))dW_{s} \\ 
 +\int_{0}^{t}\int_{H\backslash\left\{0\right\}} S(t-s)f(v,X(s))q(dv,ds)\nonumber
 \end{align*}
 holds in $[0,T]$ a.s..
\end{definition}

 \begin{definition}\label{imlddf2}
 A stochastic process $\left\{X(t),t\geq 0\right\}$ is called a strong solution of (\ref{imldeq1}) in $[0,T]$, if for all $t\leq T$\\
 (i) $X(t)$ is $\mathcal{F}_{t}$-adapted on the filtered probability space 
 ($\Omega,\mathcal{F},\left\{\mathcal{F}_{t}\right\}_{t\in [0,T]},P$),\\
 (ii) $\left\{X(t),t\geq 0\right\}$ is c\`{a}dl\`{a}g with probability one,\\
 (iii) $X(t)\in\mathcal{D}(A)$, $dt\otimes dP$ a.e., $\int_{0}^{T}\left\|AX(t)\right\|_{H}dt<\infty$  $P$ a.s.,\\
 (iv)
 \begin{align*}
 X(t)=\xi+\int_{0}^{t}(AX(s)+F(X(s)))ds+\int_{0}^{t}B(X(s))dW_{s} \\ +\int_{0}^{t}\int_{H\backslash\left\{0\right\}}f(v,X(s))q(dv,ds)
 \end{align*}
  holds in $[0,T]$ a.s..
 \end{definition}

Let $\mathcal{L}(K,H)$ be the space of all linear bounded operators from $K$ to $H$. Let $\left\{f_{j}\right\}_{j=1}^{\infty}$ be an 
 ONB in $K$ diagonailizing $Q$ and let the corresponding eigenvalues be $\left\{\lambda_{j}\right\}_{j=1}^{\infty}$.  
 Let $\mathcal{L}_{2}(K_{Q},H)$ be 
 the space of Hilbert-Schmidt operators from $K_{Q}:=Q^{1/2}K$ to $H$. (see Section 2.2 of chapter 2 of \cite{1}, Chapter 4 of 
 \cite{k17}, or \cite{P-Z}). Let $\Lambda_{2}(K_{Q},H)$ be a class of $\mathcal {L}_{2}(K_{Q},H)$-valued measurable processes, $\phi(t)$ as mapping from 
$([0,T]\times\Omega,\mathcal{B}([0,T])\otimes\mathcal{F})$ to $(\mathcal {L}_{2}(K_{Q},H),\mathcal{B}(\mathcal {L}_{2}(K_{Q},H)))$, 
adapted to the filtration $\left\{\mathcal{F}_{t}\right\}_{\left\{t\leq T\right\}}$, and satisfying the condition
$E\big[\int_{0}^{T}\left\|\phi(t)\right\|_{\mathcal {L}_{2}(K_{Q},H)}^{2}dt\big]<\infty$.\\

 Let $L^2_{T,\beta}(H)$ be the space s.t. 
 $L^2_{T,\beta}(H):=\big\{\varphi:(H\backslash\left\{0\right\})\times[0,T]\times\Omega\rightarrow H$, such that $\varphi$ is 
 jointly measurable and  $\mathcal F_{t}$-adapted for all $v\in H\backslash\{0\}$, $t\in[0,T]$ with 
$E[\int_{0}^{T}\int_{H\backslash\left\{0\right\}}\left\|\varphi(v,t)\right\|_{H}^{2}\beta(dv)dt]<\infty\big\}$.\\

Let us assume in the next Lemmma that (A1), (A2) hold.

\begin{lemma}\label{lemma1}
 a) Let $\tilde{B}(s)\in\Lambda_{2}(K_{Q},H)$ with  
$E[\int_{0}^{T}\big\|\tilde{B}(s)\big\|_{\mathcal L_{2}(K_{Q},H)}^{2}ds]<\infty$. Then for any stopping time $\tau$, there exists a constant $C_1$, 
depending on $\alpha$, $T$ s.t.
 \begin{align*}
 E\left[ \sup_{0\leq t\leq T\wedge\tau}\big\|\int_{0}^{t}S(t-s)\tilde{B}(s)dW_{s}\big\|_{H}^{2}\right]
 \leq C_{1}E \left[\int_{0}^{T\wedge\tau}\big\|\tilde{B}(s)\big\|_{\mathcal L_{2}(K_{Q},H)}^{2}ds\right].
 \end{align*}
 b) Let $\varphi\in L^2_{T,\beta}(H)$ and $\tau$ be a stopping time. Then, there exists a constant $C_2$, depending on $\alpha$, $T$ s.t.
 \begin{align*}
 E\left[\sup_{0\leq t\leq T\wedge\tau}\big\|\int_{0}^{t}\int_{H\backslash\left\{0\right\}}S(t-s)\varphi(v,s)q(dv,ds)\big\|_{H}^{2}\right] \\ 
 \leq C_{2}E\left[\int_{0}^{T\wedge\tau}\int_{H\backslash\left\{0\right\}}\left\|\varphi(v,s)\right\|_{H}^{2}\beta(dv)ds\right].
 \end{align*}
 \begin{proof}
 For the proof of the first inequality we refer Lemma 3.3(b) of \cite{1}. For the proof of the second inequality we refer Lemma 5.1.9(1) of 
 \cite{bookR}, \cite{haus.brez.zh}.

\end{proof}
\end{lemma}
 


Moreover, let us assume in this section  from this point on that (A1), (A2), (A3) holds.\\
 Let $\{\xi(t), t\in[0,T] $ be a $\mathcal F_{t}$-adapted process  for all  $t\in[0,T]$.

 Let us define 
 \begin{align}\label{I(t)}
 I(t,\xi(t))&=\int_{0}^{t}S(t-s)F(\xi(s))ds+\int_{0}^{t}S(t-s)B(\xi(s))dW_{s} \\ &
 +\int_{0}^{t}\int_{H\backslash\left\{0\right\}} S(t-s)f(v,\xi(s))q(dv,ds)\,,\quad t \in[0,T].\nonumber
 \end{align}


\begin{lemma}\label{imldlm2}
 Let $E[\sup_{0\leq s\leq T}\|\xi(s)\|^2_H]<\infty$. Then for any stopping time $\tau$
\begin{align*}
E\left[\sup_{0\leq s\leq t\wedge\tau}\left\|I(s,\xi(s))\right\|_{H}^{2}\right]\leq 
C_{3}\left(t+\int_{0}^{t}E[\sup_{0\leq u\leq s\wedge\tau}\left\|\xi(u)\right\|_{H}^{2}]ds\right),
\end{align*}
where $C_{3}$ is a constant depending on $\alpha$, $T$ and $l$. 
\begin{proof}
  For the proof we refer Lemma 3.4 of Chapter 3 of \cite{1} and Theorem 5.2.1 of \cite{bookR}. In \cite{1} it is done for the Gaussian case and in 
  \cite{bookR} it is done for the non-Gaussian case.
 \end{proof}
\end{lemma}


\begin{lemma}\label{imldlm3}
 Let $E[\sup_{0\leq s\leq T}\|\xi(s)\|^2_H]<\infty$ . Then
\begin{align*}
E\left[\sup_{0\leq s\leq t}\left\|I(s,\xi_{1}(s))-I(s,\xi_{2}(s))\right\|_{H}^{2}\right]\leq 
C_{4}\int_{0}^{t}E[\sup_{0\leq u\leq s}\left\|\xi_{1}(u)-\xi_{2}(u)\right\|_{H}^{2}]ds,
\end{align*}
where $C_{4}$ is a constant depending on $\alpha$, $T$ and $\mathcal{K}$. 
\begin{proof}
 For the proof we refer Lemma 3.5 of \cite{1} and Lemma 5.2.2 of \cite{bookR}. In \cite{1} it is done for the Gaussian case and in 
  \cite{bookR} it is done for the non-Gaussian case.
 \end{proof}
\end{lemma}


 \subsection{Existence and uniqueness of the mild solutions}\label{exs-unq}

Let $(D[0,T],H)$ be the space of c\`{a}dl\`{a}g functions defined on $[0,T]$ and with values in $H$, with the sup norm 
$\left\|.\right\|_{\infty}:=\sup_{t\in[0,T]}\left\|.\right\|_{H}$. Let $\mathcal{H}_2^T$ denote the space of $(D[0,T],H)$-valued 
random processes $\xi(t)$, which are jointly measurable, adapted to the filtration 
$\left\{\mathcal{F}_t\right\}_{t\in[0,T]}$, with $E[\sup_{0\leq s\leq T}\left\|\xi(s)\right\|_{H}^{2}]<\infty$. The space $\mathcal{H}_2^T$, 
equipped with the norm $\|X\|_{\mathcal{H}_2^T}:=\big(E[\sup_{0\leq s\leq T}\left\|X(s)\right\|_{H}^{2}]\big)^{1/2}$ is a Banach space 
(see Section 4.1 of \cite{bookR}).

\begin{theorem}\label{imldtm1}
Let $\{S(t),t\geq0\}$ be a pseudo-contraction semigroup generated by $A$ satisfying assumption (A1). Suppose assumption (A2) holds. 
Let the coefficients $F$, $B$, $f$ satisfy assumption (A3),   and assume that
$E[\left\|X(0)\right\|_{H}^{2}]<\infty$. Then for each $T>0$ equation (\ref{imldeq1}) has a unique mild solution $X^\xi\in(D[0,T],H)$ satisfying 
$E[\sup_{0\leq s\leq T}\left\|X^\xi(s)\right\|_{H}^{2}]<\infty$, i.e. the mild solution is in $\mathcal{H}_2^T$.
\begin{proof}
 For the proof we refer Theorem 3.3 of \cite{1} and Theorem 5.2.3 of \cite{bookR}. In \cite{1} it is done for the Gaussian case and in 
  \cite{bookR} it is done for the non-Gaussian case. The proofs use Lemma \ref{lemma1} and Lemma \ref{imldlm3}. In Lemma \ref{lemma1} the condition $\tilde{B}(s)\in\Lambda_{2}(K_{Q},H)$, for $\tilde{B}(s):=B(X(s))$ is a consequence of  assumption (A3). 
\end{proof}
\end{theorem}
 

\subsection{When a mild solution is a strong solution}\label{mild->strong}

\begin{theorem}\label{imldtm2}
Let $\{S(t),t\geq0\}$ be a pseudo-contraction semigroup generated by $A$ satisfying assumption (A1). Suppose assumption (A2) holds. 
Let the coefficients $F$, $B$, $f$ satisfy assumption (A3),   and assume that
$E[\left\|X(0)\right\|_{H}^{2}]<\infty$.  \\ Suppose 
 $\xi\in\mathcal{D}(A)$  a.s., $S(t-r)F(y)\in\mathcal{D}(A)$, $S(t-r)B(y)\in\mathcal{D}(A)$, $S(t-r)f(v,y)\in\mathcal{D}(A)$; for all $r<t$, $y\in H$ and $v\in H\backslash\left\{0\right\}$.\\

Then for any $T>0$ the  mild solution  $\{X(t)\}_{t\in [0,T]}$ of (\ref{imldeq1}) satisfies 
$X(t)\in\mathcal{D}(A)$  $dt\otimes dP$ a.e., and   is a strong solution.
\begin{proof} 
  Here we follow exactly the proof of Proposition 2.3 of \cite{7} and Theorem 3.2 of \cite{3}. In \cite{7}  the Gaussian case is considered and  
 in \cite{3} the non-Gaussian case is considered.\\
From the assumptions it follows

\begin{align}
\int_{0}^{t}AX(s)ds=\int_{0}^{t}AS(s)\xi ds+\int_{0}^{t}\int_{0}^{s}AS(s-r)F(X(r))drds+\int_{0}^{t}\int_{0}^{s}AS(s-r)B(X(r))dW_{r}ds\nonumber
\end{align}
\begin{align}
+\int_{0}^{t}\int_{0}^{s}\int_{H\backslash\left\{0\right\}} AS(s-r)f(v,X(r))q(dv,dr)ds,\nonumber
\end{align}

since, we know
\begin{align}
\int_{0}^{t}\int_{0}^{s}f(s-r)g(r)drds&=\int_{0}^{t}\int_{0}^{t}f(s-r)g(r)\chi_{[0,s]}(r)drds&\nonumber  \\
&=\int_{0}^{t}\int_{0}^{t}f(s-r)g(r)\chi_{[0,s]}(r)dsdr&\nonumber  \\
&=\int_{0}^{t}\int_{0}^{t}f(s-r)g(r)\chi_{[r,t]}(s)dsdr&\nonumber  \\
&=\int_{0}^{t}\int_{r}^{t}f(s-r)g(r)dsdr,&\nonumber
\end{align}
 and, by the given conditions we have
\begin{align*}
\int_{0}^{T}\int_{0}^{t}\left\|AS(t-r)F(X(r))\right\|drdt<\infty	
\end{align*}
with probability one,
\begin{align*}
\int_{0}^{T}\int_{0}^{t}\left\|AS(t-r)B(X(r))\right\|^2drdt<\infty	
\end{align*}
with probability one and
\begin{align*}
\int_{0}^{T}\int_{0}^{t}\int_{H\backslash\left\{0\right\}}E\left\|AS(t-r)f(v,X(r))\right\|^{2}\beta(dv)drdt<\infty	 
\end{align*}

with probability one. Hence, by applying Fubini theorem we get,

\begin{align}
\int_{0}^{t}AX(s)ds=\int_{0}^{t}AS(s)\xi ds+\int_{0}^{t}\int_{r}^{t}AS(s-r)F(X(r))dsdr+\int_{0}^{t}\int_{r}^{t}AS(s-r)B(X(r))dsdW_{r}\nonumber
\end{align}
\begin{align}
+\int_{0}^{t}\int_{H\backslash\left\{0\right\}}\int_{r}^{t} AS(s-r)f(v,X(r))dsq(dv,dr).\nonumber
\end{align}
(For the stochastic Fubini theorem we refer theorem (2.8) of \cite{1} and Appendix A of \cite{bookR}). Now we apply the formula,

\begin{align*}
\int_{0}^{t}AS(s)\xi ds=S(t)\xi-\xi;
\end{align*}

when $\xi\in\mathcal{D}(A)$ a.s.. Hence $AX(t)$ is integrable with probability one and
\begin{align}
\int_{0}^{t}AX(s)ds=S(t)\xi-\xi+\int_{0}^{t}S(t-r)F(X(r))dr-\int_{0}^{t}F(X(r))dr\nonumber
\end{align}
\begin{align}
+\int_{0}^{t}S(t-r)B(X(r))dW_{r}-\int_{0}^{t}B(X(r))dW_{r}\nonumber
\end{align}
\begin{align}
+\int_{0}^{t}\int_{H\backslash\left\{0\right\}}S(t-r)f(v,X(r))q(dv,dr)-\int_{0}^{t}\int_{H\backslash\left\{0\right\}}f(v,X(r))q(dv,dr).\nonumber
\end{align}
Hence
\begin{align}
\int_{0}^{t}AX(r)dr=X(t)-\xi-\int_{0}^{t}F(X(r))dr-\int_{0}^{t}B(X(r))dW_{r}\nonumber
\end{align}
\begin{align}
-\int_{0}^{t}\int_{H\backslash\left\{0\right\}}f(v,X(r))q(dv,dr).\nonumber
\end{align}
Therefore
\begin{align}
X(t)=\xi+\int_{0}^{t}AX(r)dr+\int_{0}^{t}F(X(r))dr+\int_{0}^{t}B(X(r))dW_{r}\nonumber
\end{align}
\begin{align}
+\int_{0}^{t}\int_{H\backslash\left\{0\right\}}f(v,X(r))q(dv,dr).\nonumber
\end{align}

By Definition \ref{imlddf2}, $\{X(t),t\in [0,T]\}$ is a strong solution of equation (\ref{imldeq1}).

\end{proof}
\end{theorem}

 
 
 

\subsection{Approximating a mild solution by the strong solutions}\label{yosidaSubc}
 We assume that the assumptions (A1), (A2), (A3) are satisfied for equation (\ref{imldeq1}).
 Let us  consider the approximating system of equation (\ref{imldeq1}),
 \begin{align*}
 dX(t)=(AX(t)+R_{n}F(X(t)))dt+R_{n}B(X(t))dW_{t}+\int_{H\backslash\left\{0\right\}} R_{n}f(v,X(t))q(dv,dt);
 \end{align*}
 \begin{align}\label{imldeq2}
 X(0)=\xi\in\mathcal{D}(A) \, a.s..
 \end{align}

 where $R_{n}=nR(n,A)$, and  $R(n,A)=(nI-A)^{-1}$ denotes the resolvent of $A$ evaluted at $n\in \mathbb{N}$, with   
  $n\in\rho(A)$,  $\rho(A)$ denoting the resolvent set of $A$. We have $R_{n}:H\rightarrow\mathcal{D}(A)$ and $A_{n}=AR_{n}$ are the Yosida 
 approximations of $A$ (see Chapter 1 of \cite{1}).  \\

 By applying Theorem \ref{imldtm1}, we conclude that equation (\ref{imldeq2}) has a unique mild solution, 
denoted by 
 $X_{n}^{\xi}(t)$. Then
\begin{align}
X_{n}^{\xi}(t)=S(t)\xi+\int_{0}^{t}S(t-s)R_{n}F(X_{n}^{\xi}(s))ds+\int_{0}^{t}S(t-s)R_{n}B(X_{n}^{\xi}(s))dW_{s}\nonumber
\end{align}
\begin{align}\label{imldeq3}
+\int_{0}^{t}\int_{H\backslash\left\{0\right\}} S(t-s)R_{n}f(v,X_{n}^{\xi}(s))q(dv,ds).
\end{align}

 Since the range $\mathcal{R}(R(n,A))\subset\mathcal{D}(A)$ (see Chapter 1 of \cite{1}) and the conditions of Theorem \ref{imldtm2} are  
 satisfied,  we  conclude that $X_{n}^{\xi}(t)\in\mathcal{D}(A)$ is also a strong solution.\\

Now we are in a position to approximate the mild solution of equation (\ref{imldeq1}) by the strong solutions of equation 
(\ref{imldeq2}). The mild solution of equation (\ref{imldeq1}), say $X^{\xi}(t)$, satisfies by definition 
\begin{align}
X^{\xi}(t)=S(t)\xi+\int_{0}^{t}S(t-s)F(X^{\xi}(s))ds+\int_{0}^{t}S(t-s)B(X^{\xi}(s))dW_{s}\nonumber
\end{align}
\begin{align}\label{imldeq4}
+\int_{0}^{t}\int_{H\backslash\left\{0\right\}} S(t-s)f(v,X^{\xi}(s))q(dv,ds).
\end{align}


\begin{theorem}\label{imldtm3}  
 The stochastic partial differential equation (\ref{imldeq2}) has a unique strong solution 
 $\left\{X_{n}^{\xi}(t), t\geq0\right\}$ in $D([0,T], L_{2}((\Omega,\mathcal{F}, P), H)$ for $T$ finite. Moreover 

\begin{align}\label{imldeq5}
 \lim_{n\rightarrow\infty}E\left[\sup_{0\leq t\leq T}\left\|X_{n}^{\xi}(t)-X^{\xi}(t)\right\|_{H}^{2}\right]=0,
\end{align}
where  $\left\{X^{\xi}(t),t\geq0\right\}$ is the mild solution of equation (\ref{imldeq1}).
\begin{proof}
 In Theorem \ref{imldtm1} we have already proved that there exists a unique solution of \eqref{imldeq2} in $D([0,T], L_{2}((\Omega,\mathcal{F}, P), H)$ and by 
 Theorem \ref{imldtm2} this is also a strong solution. Now we will prove \eqref{imldeq5}. We have 
\begin{align}
E[\sup_{0\leq t\leq T}\left\|X_{n}^{\xi}(t)-X^{\xi}(t)\right\|_{H}^{2}]\nonumber
\end{align}
\begin{align*}
=&E[\sup_{0\leq t\leq T}\big\|\int_{0}^{t}S(t-s)(R_{n}F(X_{n}^{\xi}(s))-F(X^{\xi}(s)))ds&\\
 &+\int_{0}^{t}S(t-s)(R_{n}B(X_{n}^{\xi}(s))-B(X^{\xi}(s)))dW_{s}&\\
 &+\int_{0}^{t}\int_{H\backslash\left\{0\right\}}S(t-s)(R_{n}f(v,X_{n}^{\xi}(s))-f(v,X^{\xi}(s)))q(dv,ds)\big\|_{H}^{2}]&
\end{align*}
 
 \begin{align*}
\leq&C\big\{E[\sup_{0\leq s\leq t}\big\|\int_{0}^{s}S(s-r)R_{n}(F(X_{n}^{\xi}(r))-F(X^{\xi}(r)))dr\big\|_{H}^{2}]&\\
 &+E[\sup_{0\leq s\leq t}\big\|\int_{0}^{s}S(s-r)(R_{n}-I)F(X^{\xi}(r))dr\big\|_{H}^{2}]&\\
 &+E[\sup_{0\leq s\leq t}\big\|\int_{0}^{s}S(s-r)R_{n}(B(X_{n}^{\xi}(r))-B(X^{\xi}(r)))dW_{r}\big\|_{H}^{2}]&\\
 &+C_{1}E[\int_{0}^{t}\big\|(R_{n}-I)B(X^{\xi}(r))\big\|_{\mathcal L_{2}(K_{Q},H)}^{2}dr]&\\
 &+E[\sup_{0\leq s\leq t}\big\|\int_{0}^{s}\int_{H\backslash\left\{0\right\}}S(s-r)R_{n}(f(v,X_{n}^{\xi}(r))-f(v,X^{\xi}(r)))q(dv,dr)\big\|_{H}^{2}]&\\
 &+C_{2}E[\int_{0}^{t}\int_{H\backslash\left\{0\right\}}\big\|(R_{n}-I)f(v,X^{\xi}(r))\big\|_{H}^{2}\beta(dv)dr]\big\},&
\end{align*}

where $C$, $C_{1}$, $C_{2}$ are constants depending on $\alpha$ and $T$. By Lemma \ref{imldlm3}, the first, third and fifth summands are bounded by 
$\mathcal{G}_{1}\mathcal{K}\int_{0}^{t}E[\sup_{0\leq r\leq s}\left\|X_{n}^{\xi}(r)-X^{\xi}(r)\right\|_{H}^{2}]dr$, for $n>n_{0}$ ($n_{0}$ sufficiently 
large), where $\mathcal{G}_{1}$ is a constant which depends on 
$\sup_{0\leq t\leq T}\left\|S(t)\right\|_{\mathcal{L}(H)}$ and $\sup_{n>n_{0}}\left\|R_{n}\right\|_{\mathcal{L}(H)}$ and $\mathcal{K}$ is the 
Lipschitz constant in (A3).   \\

By the properties of $R_{n}$, the integrands in the second, fourth and sixth summands converge to zero, as $n\to \infty$. The integrands are bounded by 
$\mathcal{G}_{2}l(1+\left\|X^{\xi}(r)\right\|_{H}^{2})$ (by condition (A3)) for some constant $\mathcal{G}_{2}$ depending on $\left\|S(t)\right\|_{\mathcal{L}(H)}$ 
and $\left\|R_{n}\right\|_{\mathcal{L}(H)}$, and the constant $l$ is the linear growth condition appearing in (A3). So by Lebesgue Dominated Convergence Theorem the 
integrals converge to zero as $n\rightarrow\infty$. Therefore there exists $\epsilon>0$ s.t. for sufficiently large $n$ each of the 
three summands are less or equal $\epsilon$. 
So for sufficiently large $n$, 
\begin{align*}
E[\sup_{0\leq t\leq T}\left\|X_{n}^{\xi}(t)-X^{\xi}(t)\right\|_{H}^{2}]\leq3\mathcal{G}_{1}\mathcal{K}
\int_{0}^{t}E[\sup_{0\leq r\leq s}\left\|X_{n}^{\xi}(r)-X^{\xi}(r)\right\|_{H}^{2}]dr+3\epsilon.
\end{align*}
By Gronwall's lemma (for sufficiently large $n$),
\begin{align*}
E[\sup_{0\leq t\leq T}\left\|X_{n}^{\xi}(t)-X^{\xi}(t)\right\|_{H}^{2}]\leq3\epsilon e^{3\mathcal{G}_{1}\mathcal{K}t}.
\end{align*}
From this  we conclude that
\begin{align*}
\lim_{n\rightarrow\infty}E[\sup_{0\leq t\leq T}\left\|X_{n}^{\xi}(t)-X^{\xi}(t)\right\|_{H}^{2}]=0.	
\end{align*}

\end{proof}
\end{theorem}

We call  $\left\{X_{n}^{\xi}(t)\right\}$ in the above theorem the Yosida approximation of the mild solution of (\ref{imldeq1}).

   \section{It\^o formula for strong and mild solutions and applications}\label{mthm}
 
 In this Section (including all subsections) we assume that the assumptions (A1), (A2), (A3) are satisfied for equation (\ref{imldeq1}).\\
 
\subsection{It\^o formula for strong solutions}
Let  all the assumption in Theorem \ref{imldtm1} be satisfied, so that a unique mild solution of (\ref{imldeq1}) exists in $\mathcal{H}_2^T$ for all $T>0$. Let us assume here that  $\xi \in \mathcal{D}(A)$ a.s. and  also that $\left\{X^\xi(t),t\geq 0\right\}$ is a strong solution of (\ref{imldeq1}). Under these assumptions  the It\^o formula known for solutions of Banach valued stochastic differential equations with Gaussian and L\'evy noise hold also for the strong solution $\left\{X^\xi(t),t\geq 0\right\}$ of (\ref{imldeq1}). We recall this result here in the more general statement obtained in \cite{9} (see Remark \ref{Rem1} for the more known statement of It\^o formula). We need to recall some notation and definition:

  Let $C^{1,2}([0,T]\times H)$ ($C_b^{1,2}([0,T]\times H)$) denote the class of real valued continuous functions $\Psi$ on $[0,T]\times H$ with continuous (uniformly bounded)  Fr\'echet derivatives $\partial_s\Psi(s,x)$, 
  $\partial_x\Psi(s,x)$, $\partial_s\partial_x\Psi(s,x)$, $\partial_x\partial_s\Psi(s,x)$ and $\partial_x\partial_x\Psi(s,x)$.
  
 \begin{definition}\label{quasi}
 We call a continuous, non-decreasing function $h:\mathbb{R}_+\rightarrow\mathbb{R}_+$ quasi-sublinear if there is a constant $C>0$ such that
 \begin{align*}
 h(x+y)\leq C(h(x)+h(y)),\ \ \ \ x,y\in\mathbb{R}_+,
 \end{align*}
 \begin{align*}
 h(xy)\leq Ch(x)h(y),\ \ \ \ x,y\in\mathbb{R}_+.
 \end{align*}
 \end{definition}

   Let $\Psi\in C^{1,2}([0,T]\times H)$. Assume moreover that the following 
   conditions hold:\\ 
  \begin{align}\label{bd-h-1}
  \left\|\Psi_x(s,x)\right\|_H\leq h_{1}(\left\|x\right\|_H)	
  \end{align}
   and
  \begin{align}\label{bd-h-2}
  \left\|\Psi_{xx}(s,x)\right\|_{\mathcal{L}(H)}\leq h_{2}(\left\|x\right\|_H). 	
  \end{align}
 
for  $h_1,h_2:\mathbb{R}_+\rightarrow\mathbb{R}_+$ quasi-sublinear functions such that
\begin{align}\label{bd-h-3}
 \int_0^T\int_{H\backslash\left\{0\right\}}\left\|f(v,X(s))\right\|^2\beta(dv)ds 
 +  \int_0^T\int_{H\backslash\left\{0\right\}}h_1(\left\|f(v,X(s))\right\|)^2\left\|f(v,X(s))\right\|^2\beta(dv)ds \\
+  \int_0^T\int_{H\backslash\left\{0\right\}}h_2(\left\|f(v,X(s))\right\|)\left\|f(v,X(s))\right\|^2\beta(dv)ds <\infty\nonumber
\end{align}
 $P$-a.s., for all $T>0$. Then the  following  usual It\^o formula holds:
 
 \begin{align}\label{24-9-15}
 &\Psi(t,X^{\xi}(t))-\Psi(0,\xi)=\int_{0}^{t}(\Psi_{s}(s,X^{\xi}(s))+\mathcal{L}\Psi(s,X^{\xi}(s)))ds&\\
 &+\int_{0}^{t}\left\langle \Psi_{x}(s,X^{\xi}(s)),B(X^{\xi}(s))dW_{s}\right\rangle_{H}&\nonumber\\
 &+\int_{0}^{t}\int_{H\backslash\left\{0\right\}}\left[\Psi(s,X^{\xi}(s)+f(v,X^{\xi}(s)))-\Psi(s,X^{\xi}(s))\right]q(dv,ds).&\nonumber
 \end{align}
 with
 \begin{align}\label{defnL}
  \mathcal{L}\Psi(s,x):=&\left\langle \Psi_{x}(s,x),Ax+F(x)\right\rangle_{H}+\frac{1}{2} tr(\Psi_{xx}(s,x)(B(x))Q(B(x))^{*})&\nonumber\\
  &+\int_{H\backslash\left\{0\right\}}\left[\Psi(s,x+f(v,x))-\Psi(s,x)-\left\langle \Psi_{x}(s,x),f(v,x)\right\rangle_{H}\right]\beta(dv)\quad s\in \mathbb{R}, x\in \mathcal{D}(A).&
  \end{align}
 This follows from the It\^o formula for strong solutions of SDEs with L\'evy noise analysed in  \cite{9} and  further in Theorem 3.7.2 in \cite{bookR} (for the more general  case of non deterministic coefficients, like e.g.  assumed here in Theorem \ref{imldtm1}).
 \begin{remark}\label{detcondition}
 A sufficient condition for \eqref{bd-h-3} to hold is 
 \begin{align}\label{bd-h-3det}
  \int_{H\backslash\left\{0\right\}}\left\|f(v,x)\right\|^2\beta(dv) 
  +  \int_{H\backslash\left\{0\right\}}(h_1(\left\|f(v,x)\right\|))^2\left\|f(v,x)\right\|^2\beta(dv) \\
 +  \int_{H\backslash\left\{0\right\}}h_2(\left\|f(v,x)\right\|)\left\|f(v,x))\right\|^2\beta(dv) <\infty \quad  \forall x\in H,\nonumber
 \end{align}
as we assume the hypothesis in Theorem \ref{imldtm1} to be satisfied, so that 
$E[\sup_{0\leq s\leq T}\left\|X^\xi(s)\right\|_{H}^{2}]<\infty$.
 \end{remark}
  \begin{remark}\label{Rem1} if $\Psi\in C_b^{1,2}([0,T]\times H)$ then the assumptions  \eqref{bd-h-1},  \eqref{bd-h-2} and  \eqref{bd-h-3} can be skipped, as we are in the case considered in \cite {10}, \cite{iki-wat}. In fact  for inequality   \eqref{bd-h-1}, \eqref{bd-h-2}, and  \eqref{bd-h-3} we can choose  the quasi linear functions $h_1$ and $h_2$ to be  constants. In particular    \eqref{bd-h-3} holds thanks to the assumption (A3) (see the first Example at page 184 in \cite{9}).
   \end{remark}

\subsection{It\^o formula for Yosida approximations of mild  solutions}

 Let    all the assumption in Theorem \ref{imldtm1} be satisfied and $\xi \in \mathcal{D}(A)$ a.s.. Let   $\left\{X_{n}^{\xi}(t)\right\}$ be the Yosida approximation of the mild solution of (\ref{imldeq1}). Assume that the quasi-linear functions satisfy \eqref{bd-h-1},   \eqref{bd-h-2} and \eqref{bd-h-3}. 
Then due to Theorem \ref{imldtm3} and  Theorem 3.7.2 in \cite{bookR} the following usual  It\^o formula holds:
 \begin{align}\label{imldeq6}
 &\Psi(t,X_{n}^{\xi}(t))-\Psi(0,\xi)=\int_{0}^{t}(\Psi_{s}(s,X_{n}^{\xi}(s))+\mathcal{L}_{n}\Psi(s,X_{n}^{\xi}(s)))ds&\nonumber\\
 &+\int_{0}^{t}\left\langle \Psi_{x}(s,X_{n}^{\xi}(s)),R_{n}B(X_{n}^{\xi}(s))dW_{s}\right\rangle_{H}&\\
 &+\int_{0}^{t}\int_{H\backslash\left\{0\right\}}\left[\Psi(s,X_{n}^{\xi}(s)+R_{n}f(v,X_{n}^{\xi}(s)))-\Psi(s,X_{n}^{\xi}(s))\right]q(dv,ds),&\nonumber
\end{align}

where

\begin{align}\label{imldeq7}
 &\mathcal{L}_{n}\Psi(s,x)=\left\langle \Psi_{x}(s,x),Ax+R_{n}F(X_{n}^{\xi}(s))\right\rangle_{H}&\\
 &+\frac{1}{2} tr(\Psi_{xx}(s,x)(R_{n}B(X_{n}^{\xi}(s)))Q(R_{n}B(X_{n}^{\xi}(s)))^{*})&\nonumber\\
 &+\int_{H\backslash\left\{0\right\}}\left[\Psi(s,x+R_{n}f(v,x)-\Psi(s,x)
   -\left\langle \Psi_{x}(s,x),R_{n}f(v,x)\right\rangle_{H}\right]\beta(dv)\quad s\in \mathbb{R}, x\in \mathcal{D}(A).&\nonumber
 \end{align}

 \begin{remark}\label{Rem2} If $\Psi\in C_b^{1,2}([0,T]\times H)$, then  \eqref{bd-h-1},  \eqref{bd-h-2} and  \eqref{bd-h-3} can be ignored,   the It\^o formula \eqref{imldeq6} still holding (See Remark \ref{Rem1}) 
    \end{remark}






 \subsection{An It\^o formula for mild solutions  obtained through Yosida approximation}\label{111sub}
  


 By using Yosida approximations we will derive an It\^o formula for mild solutions of \eqref{imldeq1} which are not necessarily strong solutions.
\begin{theorem}\label{imldtm4}
Let $\{S(t),t\geq0\}$ be a pseudo-contraction semigroup generated by $A$ satisfying assumption (A1). Suppose moreover that assumption (A2) holds. 
Let the coefficients $F$, $B$, $f$ satisfy assumption (A3),   and assume that $\xi \in \mathcal{D}(A)$ a.s..  Moreover assume that \eqref{bd-h-1}, \eqref{bd-h-2},  \eqref{bd-h-3} hold.

 Then the following It\^o Formula for the mild solution $X^\xi$ of  equation (\ref{imldeq1})  holds $P$-a.s. for all $t\in[0,T]$  
\begin{align}\label{ito-lim-forml}
 &\lim_{n\rightarrow\infty}\int_{0}^{t}\left\langle\Psi_{x}(s,X_{n}^{\xi}(s)),AX_{n}^{\xi}(s)\right\rangle_{H}ds&\\
 &=\Psi(t,X^{\xi}(t))-\Psi(0,\xi)\nonumber-\int_{0}^{t}(\Psi_{s}(s,X^{\xi}(s)))ds
 -\int_{0}^{t}\left\langle\Psi_{x}(s,X^{\xi}(s)),F(X^{\xi}(s))\right\rangle_{H}ds&\nonumber\\
 &-\int_{0}^{t}\frac{1}{2} tr(\Psi_{xx}(s,X^{\xi}(s))(B(X^{\xi}(s)))Q(B(X^{\xi}(s)))^{*})ds&\nonumber\\
 &-\int_{0}^{t}\int_{H\backslash\left\{0\right\}}\left[\Psi(s,X^{\xi}(s)+f(v,X^{\xi}(s)))-\Psi(s,X^{\xi}(s))
  -\left\langle \Psi_{x}(s,X^{\xi}(s)),f(v,X^{\xi}(s))\right\rangle_{H}\right]\beta(dv)ds&\nonumber\\
 &-\int_{0}^{t}\left\langle \Psi_{x}(s,X^{\xi}(s)),B(X^{\xi}(s))dW_{s}\right\rangle_{H}&\nonumber\\
 &-\int_{0}^{t}\int_{H\backslash\left\{0\right\}}\left[\Psi(s,X^{\xi}(s)+f(v,X^{\xi}(s)))-\Psi(s,X^{\xi}(s))\right]q(dv,ds).&\nonumber
\end{align}
\end{theorem}

 We will use in the proof of Theorem \ref{imldtm4} below the It\^o formula for the Yosida approximations of  the mild solution $X^\xi$ of  equation (\ref{imldeq1}). It is then easy checked that due Remark \ref{Rem2} the following Corollary is a direct consequence of the proof of Theorem \ref{imldtm4}. 

\begin{corollary}\label{Rem3}
Let $\{S(t),t\geq0\}$ be a pseudo-contraction semigroup generated by $A$ satisfying assumption (A1). Suppose assumption (A2) holds. 
Let the coefficients $F$, $B$, $f$ satisfy assumption (A3),   and assume that $\xi \in \mathcal{D}(A)$ a.s.. Then for each  $\Psi\in C_b^{1,2}([0,T]\times H)$ the 
It\^o Formula \eqref{ito-lim-forml} holds
\end{corollary}
Here we will prove Theorem \ref{imldtm4}:
\begin{proof}
 From \eqref{imldeq6},
 \begin{align*}
 &\Psi(t,X_{n}^{\xi}(t))-\Psi(0,\xi)&\\
  &=\int_{0}^{t}(\Psi_{s}(s,X_{n}^{\xi}(s))+\mathcal{L}_{n}\Psi(s,X_{n}^{\xi}(s)))ds
    +\int_{0}^{t}\left\langle \Psi_{x}(s,X_{n}^{\xi}(s)),R_{n}B(X_{n}^{\xi}(s))dW_{s}\right\rangle_{H}&\\
 &+\int_{0}^{t}\int_{H\backslash\left\{0\right\}}\left[\Psi(s,X_{n}^{\xi}(s)+R_{n}f(x,X_{n}^{\xi}(s)))
   -\Psi(s,X_{n}^{\xi}(s))\right]q(dx,ds).&
 \end{align*}


Now we substite  $\mathcal{L}_{n}(\Psi(s,X_{n}^{\xi}(s)))$,
\begin{align}\label{imldeq8}
=&\int_{0}^{t}(\Psi_{s}(s,X_{n}^{\xi}(s)))ds
   +\int_{0}^{t}\left\langle \Psi_{x}(s,X_{n}^{\xi}(s)),AX_{n}^{\xi}(s)+R_{n}F(X_{n}^{\xi}(s))\right\rangle_{H}ds&\\
 &+\int_{0}^{t}\frac{1}{2} tr(\Psi_{xx}(s,X_{n}^{\xi}(s))(R_{n}B(X_{n}^{\xi}(s)))Q(R_{n}B(X_{n}^{\xi}(s)))^{*})ds&\nonumber\\
 &+\int_{0}^{t}\int_{H\backslash\left\{0\right\}}\left[\Psi(s,X_{n}^{\xi}(s)+R_{n}f(v,X_{n}^{\xi}(s)))-\Psi(s,X_{n}^{\xi}(s))
    -\left\langle \Psi_{x}(s,X_{n}^{\xi}(s)),R_{n}f(v,X_{n}^{\xi}(s))\right\rangle_{H}\right]\beta(dv)ds&\nonumber\\
 &+\int_{0}^{t}\left\langle \Psi_{x}(s,X_{n}^{\xi}(s)),R_{n}B(X_{n}^{\xi}(s))dW_{s}\right\rangle_{H}&\nonumber\\
 &+\int_{0}^{t}\int_{H\backslash\left\{0\right\}}\left[\Psi(s,X_{n}^{\xi}(s)+R_{n}f(v,X_{n}^{\xi}(s)))-\Psi(s,X_{n}^{\xi}(s))\right]q(dv,ds).&\nonumber
\end{align}



Now our task is to show that the above equation converges $P$-a.s.(term by term) and also to find the limit.\\



The convergence in Theorem \ref{imldtm3} (equation (\ref{imldeq5})) allows us to choose a subsequence $X_{n_{k}}^{\xi}$ such that,

\begin{center}
$X_{n_{k}}^{\xi}(t)\rightarrow X^{\xi}(t)$, $0\leq t\leq T$, $P$-a.s.
\end{center}
We will denote such a subsequence  again by $X_{n}^{\xi}$.

In fact, we can say that
\begin{align}\label{sup}
\sup_{0\leq t\leq T}\left\| X_n(t)-X(t)\right\|_H\to 0, 
\end{align}
$P$ a.s.. This implies that the set
\begin{align}\label{SetA}
S=\left\{ X_n(t),\; X(t):\; n=1,2...,\, 0\leq t\leq T\right\}
\end{align}
is bounded in $H$, hence all the values of $\Psi$ and its derivatives evaluated on $S$ are bounded by some constant. 
Now we are ready to show the term by term convergence of equation (\ref{imldeq8}).\\

First consider the first term of the L.H.S. of eq. (\ref{imldeq8}). Since $\Psi$ is continuous, from (\ref{sup}) we conclude that
\begin{align*}
\Psi(t,X_{n}^{\xi}(t))\rightarrow\Psi(t,X^{\xi}(t)),
\end{align*}
$P$-a.s.\\

Now consider the first term of the R.H.S. of eq. (\ref{imldeq8}). $\Psi_{s}$ is continuous, $\Psi_{s}(s,X_{n}^{\xi}(s))<C$ 
by equation (\ref{sup}). So by applying Lebesgue dominated convergence Theorem we conclude that

 


\begin{align*}
\int_{0}^{t}(\Psi_{s}(s,X_{n}^{\xi}(s)))ds\rightarrow\int_{0}^{t}(\Psi_{s}(s,X^{\xi}(s)))ds,
\end{align*}
$P$-a.s.\\

Now consider the second term of the R.H.S. of eq. \eqref{imldeq8},

\begin{align*}
 &\int_{0}^{t}\left\langle \Psi_{x}(s,X_{n}^{\xi}(s)),AX_{n}^{\xi}(s)+R_{n}F(X_{n}^{\xi}(s))\right\rangle_{H}ds&\\
 &=\int_{0}^{t}\left\langle \Psi_{x}(s,X_{n}^{\xi}(s)),AX_{n}^{\xi}(s)\right\rangle_{H}ds
 +\int_{0}^{t}\left\langle \Psi_{x}(s,X_{n}^{\xi}(s)),R_{n}F(X_{n}^{\xi}(s))\right\rangle_{H}ds.&	
\end{align*}

Since $\Psi_{x}$ is continuous, by (\ref{sup}) we get $\Psi_{x}(s,X_{n}^{\xi}(s))\rightarrow\Psi_{x}(s,X^{\xi}(s))$. 
Since $F$ is continuous and $R_n(F(X_n^\xi(s))$ is a double sequence, therefore we have


 \begin{eqnarray}\label{DoubleSeqArg}
\left\| R_n\left( F\left( X_n(s)\right)\right) - F(X(s))\right\|_H&\leq &\left\| R_n\left( F\left( X_n(s)\right) - F(X(s))\right)\right\|_H
+\left\| R_n\left( F\left( X(s)\right)\right) - F(X(s))\right\|_H\nonumber\\
&\leq& \left\| R_n\right\|_H\left\| F\left( X_n(s)\right) - F(X(s))\right\|_H+\left\|\left( R_n-I\right)F(X(s))\right\|_H.
\end{eqnarray}
Therefore $R_n\left( F\left( X_n(s)\right)\right) \to F(X(s))$ because of the uniform boundedness of $\left\|R_{n}\right\|_{\mathcal {L}(H)}$, 
and the convergence of $(R_{n}-I)x\rightarrow 0$. So, by \eqref{sup} and Lebesgue dominated convergence theorem,



\begin{align}
\int_{0}^{t}\left\langle \Psi_{x}(s,X_{n}^{\xi}(s)),R_{n}F(X_{n}^{\xi}(s))\right\rangle_{H}ds\rightarrow\int_{0}^{t}\left\langle \Psi_{x}(s,X^{\xi}(s)),F(X^{\xi}(s))\right\rangle_{H}ds\nonumber	
\end{align}

$P$-a.s..  \\

We will discuss the convergence of the term,
\begin{align*}
\int_{0}^{t}\left\langle \Psi_{x}(s,X_{n}^{\xi}(s)),AX_{n}^{\xi}(s)\right\rangle_{H}ds	
\end{align*}
at the end.\\

Now consider the third term of the R.H.S. of eq. (\ref{imldeq8}),
\begin{align*}
\int_{0}^{t}\frac{1}{2} tr(\Psi_{xx}(s,X_{n}^{\xi}(s))(R_{n}B(X_{n}^{\xi}(s)))Q(R_{n}B(X_{n}^{\xi}(s)))^{*})ds.
\end{align*}

We have

\begin{align*}
 &tr(\Psi_{xx}(s,X_{n}^{\xi}(s))(R_{n}B(X_{n}^{\xi}(s)))Q(R_{n}B(X_{n}^{\xi}(s)))^{*})&\\
 &=tr((R_{n}B(X_{n}^{\xi}(s)))^{*}\Psi_{xx}(s,X_{n}^{\xi}(s))(R_{n}B(X_{n}^{\xi}(s)))Q)&\\
 &=\sum_{j=1}^{\infty}\lambda_{j}\left\langle \Psi_{xx}(s,X_{n}^{\xi}(s))(R_{n}B(X_{n}^{\xi}(s)))f_{j},(R_{n}B(X_{n}^{\xi}(s)))f_{j}
   \right\rangle_{H}.&
\end{align*}

Here we used the property that, for a symmetric operator $T\in\mathcal {L}(H)$ and\\ $\phi\in\mathcal {L}(K,H)$,

\begin{align*}
tr(T\phi Q\phi^{*})=tr(\phi^{*}T\phi Q).
\end{align*}

 $\Psi_{xx}$ being continuous, $B$  is continuous, $\left\|R_{n}\right\|_{\mathcal {L}(H)}$ is uniformly bounded and having the convergence 
 of $(R_{n}-I)x\rightarrow 0$, by a similar calculation as in~\eqref{DoubleSeqArg} we  deduce that

 \begin{align*}
 &\left\langle \Psi_{xx}(s,X_{n}^{\xi}(s))(R_{n}B(X_{n}^{\xi}(s)))f_{j},(R_{n}B(X_{n}^{\xi}(s)))f_{j}\right\rangle_{H}&\\
 &\rightarrow\left\langle \Psi_{xx}(s,X^{\xi}(s))(B(X^{\xi}(s)))f_{j},(B(X^{\xi}(s)))f_{j}\right\rangle_{H}.&
 \end{align*}

Hence,
 \begin{align*}
 &tr(\Psi_{xx}(s,X_{n}^{\xi}(s))(R_{n}B(X_{n}^{\xi}(s)))Q(R_{n}B(X_{n}^{\xi}(s)))^{*})&\\
 &\rightarrow tr(\Psi_{xx}(s,X^{\xi}(s))(B(X^{\xi}(s)))Q(B(X^{\xi}(s)))^{*}).&
 \end{align*}
Also we have,
\begin{align*}
 tr(\Psi_{xx}(s,X_{n}^{\xi}(s))(R_{n}B(X_{n}^{\xi}(s)))Q(R_{n}B(X_{n}^{\xi}(s)))^{*})
 &\leq\|\Psi_{xx}(s,X_{n}^{\xi}(s))\|\|R_{n}B(X_{n}^{\xi}(s))\|^2&\\
 \text{by (A3)}\ \ &\leq\|\Psi_{xx}(s,X_{n}^{\xi}(s))\|\|R_{n}\|^2l(1+\|X_n^{\xi}(s)\|^2).&
\end{align*}

So by (\ref{sup}) and Lebesgue dominated convergence theorem we  conclude that,

 \begin{align*}
 &\int_{0}^{t}tr(\Psi_{xx}(s,X_{n}^{\xi}(s))(R_{n}B(X_{n}^{\xi}(s)))Q(R_{n}B(X_{n}^{\xi}(s)))^{*})ds&\\
 &\rightarrow \int_{0}^{t}tr(\Psi_{xx}(s,X^{\xi}(s))(B(X^{\xi}(s)))Q(B(X^{\xi}(s)))^{*})ds,&
 \end{align*}
$P$-a.s..\\





Now consider the fourth term of the R.H.S. of eq. (\ref{imldeq8}),

 \begin{align*}
 \int_{H\backslash\left\{0\right\}}\left[\Psi(s,X_{n}^{\xi}(s)+R_{n}f(v,X_{n}^{\xi}(s)))-\Psi(s,X_{n}^{\xi}(s))-\left\langle \Psi_{x}(s,X_{n}^{\xi}(s)),R_{n}f(v,X_{n}^{\xi}(s))\right\rangle_{H}\right]\beta(dv). 
\end{align*}

Using  Theorem \ref{imldtm3}, \eqref{sup}, the continuity of $\Psi$, $\Psi_x$, $f$ and $(R_n-I)x\rightarrow0$, we  conclude

\begin{align*}
  \left[\Psi(s,X_{n}^{\xi}(s)+R_{n}f(v,X_{n}^{\xi}(s)))-\Psi(s,X_{n}^{\xi}(s))-\left\langle \Psi_{x}(s,X_{n}^{\xi}(s)),R_{n}f(v,X_{n}^{\xi}(s))\right\rangle_{H}\right] 
\end{align*}
converges to
\begin{align*}
  \left[\Psi(s,X^{\xi}(s)+f(v,X^{\xi}(s)))-\Psi(s,X^{\xi}(s))-\left\langle \Psi_{x}(s,X^{\xi}(s)),f(v,X^{\xi}(s))\right\rangle_{H}\right] 
\end{align*}

$P$-a.s.. Again by Taylor's theorem, the Cauchy Schwarz inequality and assumption \eqref{bd-h-2}, we get
 
\begin{align*}
 &\int_{H\backslash\left\{0\right\}}|\Psi(s,X_{n}^{\xi}(s)+R_{n}f(v,X_{n}^{\xi}(s)))-\Psi(s,X_{n}^{\xi}(s))-\left\langle 
      \Psi_{x}(s,X_{n}^{\xi}(s)),R_{n}f(v,X_{n}^{\xi}(s))\right\rangle_{H}|\beta(dv)&\\
 &=\int_{H\backslash\left\{0\right\}}|\int_0^1\Psi_{xx}(s,X_{n}^{\xi}(s)+\theta R_{n}f(v,X_{n}^{\xi}(s)))
     \left\langle R_{n}f(v,X_{n}^{\xi}(s)),R_{n}f(v,X_{n}^{\xi}(s))\right\rangle d\theta|\beta(dv)&\\
 &\leq \int_{H\backslash\left\{0\right\}}\int_0^1\left\|\Psi_{xx}(s,X_{n}^{\xi}(s)
      +\theta R_{n}f(v,X_{n}^{\xi}(s)))\right\|\left\|R_{n}f(v,X_{n}^{\xi}(s))\right\|^2d\theta\beta(dv)&\\
 &\leq \int_{H\backslash\left\{0\right\}}\int_0^1h_2\left(\left\|X_{n}^{\xi}(s)
      +\theta R_{n}f(v,X_{n}^{\xi}(s))\right\|\right)\left\|R_{n}f(v,X_{n}^{\xi}(s))\right\|^2d\theta\beta(dv)&\\ 
 &\leq C \int_{H\backslash\left\{0\right\}}\int_0^1\left(h_2(\left\|X_{n}^{\xi}(s)\right\|)
      +Ch_2(\theta)h_2(\left\|R_{n}f(v,X_{n}^{\xi}(s))\right\|)\right)\left\|R_{n}f(v,X_{n}^{\xi}(s))\right\|^2d\theta\beta(dv)&\\
 &\leq C\int_{H\backslash\left\{0\right\}}h_2(\left\|X_{n}^{\xi}(s)\right\|)\left\|R_{n}f(v,X_{n}^{\xi}(s))\right\|^2\beta(dv)&\\
&\ \ \ \ \ \ \ \ +C^2h_2(1) \int_{H\backslash\left\{0\right\}}h_2(\left\|R_{n}f(v,X_{n}^{\xi}(s))\right\|)\left\|R_{n}f(v,X_{n}^{\xi}(s))\right\|^2\beta(dv)<\infty,& 
\end{align*}

 $P$-a.s. by condition \eqref{bd-h-3}. Since $\left\|R_n\right\|_{\mathcal{L}(H)}$ is uniformly bounded, 
therefore by Lebesgue dominated convergence theorem 

\begin{align*}
 \int_0^t\int_{H\backslash\left\{0\right\}}\left[\Psi(s,X_{n}^{\xi}(s)+R_{n}f(v,X_{n}^{\xi}(s)))-\Psi(s,X_{n}^{\xi}(s))-\left\langle \Psi_{x}(s,X_{n}^{\xi}(s)),R_{n}f(v,X_{n}^{\xi}(s))\right\rangle_{H}\right]\beta(dv)ds 
\end{align*}

converges to

\begin{align*}
\int_{0}^{t}\int_{H\backslash\left\{0\right\}}\left[\Psi(s,X^{\xi}(s)+f(v,X^{\xi}(s)))-\Psi(s,X^{\xi}(s))-\left\langle \Psi_{x}(s,X^{\xi}(s)),f(v,X^{\xi}(s))\right\rangle_{H}\right]\beta(dv)ds
\end{align*}

$P$-a.s..\\


Now consider the fifth term of the R.H.S. of eq. (\ref{imldeq8}),
\begin{align*}
\int_{0}^{t}\left\langle \Psi_{x}(s,X_{n}^{\xi}(s)),R_{n}B(X_{n}^{\xi}(s))dW_{s}\right\rangle_{H}.
\end{align*}

\begin{align*}
 &E|\int_{0}^{t}\left\langle \Psi_{x}(s,X_{n}^{\xi}(s)),R_{n}B(X_{n}^{\xi}(s))dW_{s}\right\rangle_{H}-\int_{0}^{t}
      \left\langle  \Psi_{x}(s,X^{\xi}(s)),B(X^{\xi}(s))dW_{s}\right\rangle_{H}|^{2}&\\	
 &\leq C\int_{0}^{t}E\big\|(B(X^{\xi}(s)))^{*}(\Psi_{x}(s,X_{n}^{\xi}(s)-\Psi_{x}(s,X^{\xi}(s)))\big\|_{\mathcal {L}_{2}(K_{Q},H)}^{2}ds&\\
 &\ \ \ +C\int_{0}^{t}E\big\|((B(X^{\xi}(s)))^{*}-(R_{n}B(X_{n}^{\xi}(s)))^{*})\Psi_{x}(s,X_{n}^{\xi}(s))\big\|_{\mathcal {L}_{2}(K_{Q},H)}^{2}ds&\\	
 &\leq C\int_{0}^{t}E(\big\|(B(X^{\xi}(s))\big\|_{\mathcal {L}_{2}(K_{Q},H)}^{2}\big\|\Psi_{x}(s,X_{n}^{\xi}(s))
     -\Psi_{x}(s,X^{\xi}(s))\big\|_{H}^{2})ds&\\
 &\ \ \ +C\int_{0}^{t}E(\big\|(B(X^{\xi}(s)))^{*}-(R_{n}B(X_{n}^{\xi}(s)))^{*}\big\|_{\mathcal {L}_{2}(K_{Q},H)}^{2}\big\|\Psi_{x}(s,X_{n}^{\xi}(s))
\big\|_{H}^{2})ds.&
\end{align*}

Here, the first integral converges to zero, since the first factor is an integrable process, and the second factor converges to zero almost surely, 
so we can apply Lebesgue dominated convergence theorem. The second integral is bounded by 
$M\big\|(B(X^{\xi}(s)))^{*}-(R_{n}B(X_{n}^{\xi}(s)))^{*}\big\|_{\Lambda_{2}(K_{Q},H)}^{2}$ for some constant $M$ 
(from \eqref{sup}, $\Psi_x$ is bounded by some constant), since 
$R_{n}B(X_{n}^{\xi}(s))\rightarrow B(X^{\xi}(s))$ in the space $\Lambda_{2}(K_{Q},H)$, so the second integral also converges to zero by Lebesgue dominated convergence theorem. 
Hence we  conclude that,
\begin{align*}
\int_{0}^{t}\left\langle \Psi_{x}(s,X_{n}^{\xi}(s)),R_{n}B(X_{n}^{\xi}(s))dW_{s}\right\rangle_{H}\rightarrow\int_{0}^{t}\left\langle \Psi_{x}(s,X^{\xi}(s)),B(X^{\xi}(s))dW_{s}\right\rangle_{H}
\end{align*}
 in mean square, therefore in probability.\\





Now consider the sixth term of the R.H.S. of eq. (\ref{imldeq8}), 

 \begin{align*}
\int_{0}^{t}\int_{H\backslash\left\{0\right\}}\left[\Psi(s,X_{n}^{\xi}(s)+R_{n}f(v,X_{n}^{\xi}(s)))-\Psi(s,X_{n}^{\xi}(s))\right]q(dv,ds).
\end{align*}

\begin{align*}
 &|\left\{\Psi(s,X_{n}^{\xi}(s)+R_{n}f(v,X_{n}^{\xi}(s)))-\Psi(s,X_{n}^{\xi}(s))\right\}-\left\{\Psi(s,X^{\xi}(s)+f(v,X^{\xi}(s)))
    -\Psi(s,X^{\xi}(s))\right\}|^2&\\
 &=|\left[\Psi(s,X_{n}^{\xi}(s)+R_{n}f(v,X_{n}^{\xi}(s)))-\Psi(s,X^{\xi}(s)+f(v,X^{\xi}(s)))\right]+\left[\Psi(s,X^{\xi}(s))
       -\Psi(s,X_{n}^{\xi}(s))\right]|^2&\\
 &\leq2|\Psi(s,X_{n}^{\xi}(s)+R_{n}f(v,X_{n}^{\xi}(s)))-\Psi(s,X^{\xi}(s)+f(v,X^{\xi}(s)))|^2&\\
 &\ \ \ \ \ \ \ \ \ \ \ \ +2|\Psi(s,X^{\xi}(s))-\Psi(s,X_{n}^{\xi}(s))|^2&\\
 &\leq2\left\|X_{n}^{\xi}(s)+R_{n}f(v,X_{n}^{\xi}(s))-\left\{X^{\xi}(s)+f(v,X^{\xi}(s))\right\}\right\|^2\sup_{0<\theta\leq 1}
    \left\|\Psi_x(s,\eta_1(\theta))\right\|^2&\\
 &\ \ \ \ \ \ \ \ \ \ \ \ +2\left\|X^{\xi}(s)-X_{n}^{\xi}(s)\right\|^2\sup_{0<\theta\leq 1}\left\|\Psi_x(s,\eta_2(\theta))\right\|^2&
\end{align*}

to obtain the above inequality, we used the following inequality 
\begin{align*}
 \|\Psi(x)-\Psi(y)\|\leq\|x-y\|\sup_{0<\theta\leq 1}\|\Psi_x(y+\theta(x-y))\|.
\end{align*}
Where 
\begin{align*}
 \eta_1(\theta)=X^{\xi}(s)+f(v,X^{\xi}(s))+\theta\left(X_{n}^{\xi}(s)-X^{\xi}(s)+R_{n}f(v,X_{n}^{\xi}(s))-f(v,X^{\xi}(s))\right)
\end{align*}
and
\begin{align*}
 \eta_2(\theta)=X_n^{\xi}(s)+\theta\left(X^{\xi}(s)-X_n^{\xi}(s)\right).
\end{align*}
Therefore, by using condition \eqref{bd-h-1}, we can write the above inequality is
\begin{align*}
 &|\left\{\Psi(s,X_{n}^{\xi}(s)+R_{n}f(v,X_{n}^{\xi}(s)))-\Psi(s,X_{n}^{\xi}(s))\right\}-\left\{\Psi(s,X^{\xi}(s)+f(v,X^{\xi}(s)))
    -\Psi(s,X^{\xi}(s))\right\}|^2\\
\leq & 4\left\{\left\|X_{n}^{\xi}(s)-X^{\xi}(s)\right\|^2+\left\|R_{n}f(v,X_{n}^{\xi}(s))-f(v,X^{\xi}(s))\right\|^2\right\}\sup_{0<\theta\leq 1}\{h_1(\|\eta_1(\theta)\|)\}^2\\
&+2\left\|X^{\xi}(s)-X_{n}^{\xi}(s)\right\|^2\sup_{0<\theta\leq 1}\{h_1(\left\|\eta_2(\theta)\right\|)\}^2.
\end{align*}

Now as $n\rightarrow\infty$, the R.H.S. of the above inequality converges to $0$ $P$-a.s.. Therefore

\begin{align}\label{ajge1}
\lim_{n\rightarrow\infty}\int_{0}^{t}\int_{H\backslash\left\{0\right\}}\|\left\{\Psi(s,X_{n}^{\xi}(s)+R_{n}f(v,X_{n}^{\xi}(s)))-\Psi(s,X_{n}^{\xi}(s))\right\}\\
-\left\{\Psi(s,X^{\xi}(s)+f(v,X^{\xi}(s)))-\Psi(s,X^{\xi}(s))\right\}\|^2\beta(dv)ds=0\nonumber
\end{align}

 $P$-a.s.. Again by Taylor's theorem, the Cauchy -Schwarz inequality and assumption \eqref{bd-h-1} we get

\begin{align}\label{ajge2}
&\int_{0}^{t}\int_{H\backslash\left\{0\right\}}|\Psi(s,X_{n}^{\xi}(s)+R_{n}f(v,X_{n}^{\xi}(s)))
    -\Psi(s,X_{n}^{\xi}(s))|^2\beta(dv)ds&\nonumber\\
 &=\int_{0}^{t}\int_{H\backslash\left\{0\right\}}|\int_0^1\Psi_x(s,X_{n}^{\xi}(s)+\theta R_{n}f(v,X_{n}^{\xi}(s)))
      R_{n}f(v,X_{n}^{\xi}(s))d\theta|^2\beta(dv)ds&\nonumber\\
 &\leq\int_{0}^{t}\int_{H\backslash\left\{0\right\}}\int_0^1\left\|\Psi_x(s,X_{n}^{\xi}(s)
      +\theta R_{n}f(v,X_{n}^{\xi}(s)))\right\|^2\left\|R_{n}f(v,X_{n}^{\xi}(s))\right\|^2d\theta\beta(dv)ds&\nonumber\\
 &\leq\int_{0}^{t}\int_{H\backslash\left\{0\right\}}\int_0^1h_1(\left\|X_{n}^{\xi}(s)+\theta R_{n}f(v,X_{n}^{\xi}(s))\right\|)^2
       \left\|R_{n}f(v,X_{n}^{\xi}(s))\right\|^2d\theta\beta(dv)ds&\nonumber\\
 &\leq C^2\int_{0}^{t}\int_{H\backslash\left\{0\right\}}\int_0^1\left\{h_1(\left\|X_{n}^{\xi}(s)\right\|)
    +Ch_1(\theta)h_1(\left\|R_{n}f(v,X_{n}^{\xi}(s))\right\|)\right\}^2\left\|R_{n}f(v,X_{n}^{\xi}(s))\right\|^2d\theta\beta(dv)ds&\nonumber\\
 &\leq2C^2\int_{0}^{t}\int_{H\backslash\left\{0\right\}}h_1(\left\|X_{n}^{\xi}(s)\right\|)^2\left\|R_{n}f(v,X_{n}^{\xi}(s))\right\|^2\beta(dv)ds&\nonumber\\
&\ \ \ \ +2C^4h_1(1)\int_{0}^{t}\int_{H\backslash\left\{0\right\}}h_1(\left\|R_{n}f(v,X_{n}^{\xi}(s))\right\|)^2\left\|R_{n}f(v,X_{n}^{\xi}(s))\right\|^2
\beta(dv)ds<\infty,&
\end{align}

 $P$-a.s. by the condition \eqref{bd-h-3}. Therefore from (\ref{ajge1}) and (\ref{ajge2}) we conclude that

\begin{align*}
 &\lim_{n\rightarrow\infty}\int_{0}^{t}\int_{H\backslash\left\{0\right\}}\left[\Psi(s,X_{n}^{\xi}(s)+R_{n}f(v,X_{n}^{\xi}(s)))
      -\Psi(s,X_{n}^{\xi}(s))\right]q(dv,ds)&\\
 &=\int_{0}^{t}\int_{H\backslash\left\{0\right\}}\left[\Psi(s,X^{\xi}(s)+f(v,X^{\xi}(s)))-\Psi(s,X^{\xi}(s))\right]q(dv,ds)&
\end{align*}

in probability.\\

Thus we have showed the term by term convergence of left- and right- hand sides of eq. (\ref{imldeq8}) except for the term 
$\int_{0}^{t}\left\langle\Psi_{x}(s,X_{n}^{\xi}(s)),AX_{n}^{\xi}(s)\right\rangle_{H}ds$. Now since all the terms of the eq. (\ref{imldeq8}) converge, 
so the term $\int_{0}^{t}\left\langle\Psi_{x}(s,X_{n}^{\xi}(s)),AX_{n}^{\xi}(s)\right\rangle_{H}ds$ has to converge. Where the nonstochastic integrals 
converge in $P$-a.s. sense and stochastic integrals converge in probability. In conclusion, possibly for a subsequence of left- and right- hand sides 
of eq. (\ref{imldeq8}) converges in $P$-a.s. sense for all $t\in[0,T]$. Hence we  conclude that eq. (\ref{imldeq6}) converges $P$-a.s. and 
we can write our It\^o formula for mild solutions as 
$\lim_{n\rightarrow\infty}\int_{0}^{t}\left\langle\Psi_{x}(s,X_{n}^{\xi}(s)),AX_{n}^{\xi}(s)\right\rangle_{H}ds$ i.e. \eqref{ito-lim-forml}. 
This completes the proof.

\end{proof}

\begin{remark} \label{Rem4} When the mild solution satisfies $X^{\xi}(t)\in \mathcal{D}(A)$, then $AX^{\xi}(t)$ is well defined. Hence, 
we get back the same It\^o formula for strong solutions as in  
eq. \eqref{imldeq6}, as by Theorem \ref{imldtm2} when $X^{\xi}(t)\in \mathcal{D}(A)$, then $X^{\xi}(t)$is also a strong solution.
\end{remark} 

In the next section we will use the following

\begin{corollary}\label{itocorollary}
 Assume that all assumptions in Theorem \ref{imldtm4} are satisfied. 
 Let $X_n^{\xi}$ be the Yosida approximation of the mild solution $X^{\xi}$ of \eqref{imldeq1}, i.e  let $X_n^{\xi}$ be the strong solution of \eqref{imldeq2}, then  
 \begin{equation}
 \label{diffL}
 \lim_{n\to \infty}|\mathcal{L}\Psi(s,X_n^{\xi}(s))-\mathcal{L}_n\Psi(s,X_n^{\xi}(s))|= 0 \quad  P-a.s.
 \end{equation}

\begin{proof} 
 As we assume $\xi \in \mathcal{D}(A)$ a.s.,  and make the assumptions  that \eqref{bd-h-1}, \eqref{bd-h-2},  \eqref{bd-h-3} hold,  the operators $\mathcal{L}$ in \eqref{defnL}, and  $\mathcal{L}_n$ in \eqref{imldeq7} are well defined in all their terms.
\begin{align}
 & \lim_{n\to \infty}| \mathcal{L}\Psi(s,X_n^{\xi}(s))-\mathcal{L}_n\Psi(s,X_n^{\xi}(s))| \label{limL}\\ &
 =\lim_{n\rightarrow\infty} |\left\langle \Psi_{x}(s,X_{n}^{\xi}(s)),AX_{n}^{\xi}(s)\right\rangle
   \notag\\&
 +\left\langle \Psi_{x}(s,X_{n}^{\xi}(s)),F(X^{\xi}(s))\right\rangle_{H}+ \frac{1}{2} tr(\Psi_{xx}(s,X_{n}^{s,\xi}(s))(B(X_{n}^{\xi}(s)))Q(B(X_{n}^{\xi}(s)))^{*})
  \notag \\&
 +\int_{H\backslash\left\{0\right\}}\left[\Psi(s,X_{n}^{\xi}(s)+f(v,X_{n}^{\xi}(s)))-\Psi(s,X_{n}^{\xi}(s))-\left\langle \Psi_{x}(s,X_{n}^{\xi}(s)),f(v,X_{n}^{\xi}(s))\right\rangle_{H}\right]\beta(dv)
  \notag \\&
 -\left\langle \Psi_{x}(s,X_{n}^{\xi}(s)),AX_{n}^{\xi}(s)\right\rangle-\left\langle \Psi_{x}(s,X_{n}^{\xi}(s)),R_{n}F(X_{n}^{\xi}(s))\right\rangle_{H}
  \notag \\&
 -\frac{1}{2} tr(\Psi_{xx}(s, X_{n}^{\xi}(s))(R_{n}B(X_{n}^{\xi}(s)))Q(R_{n}B(X_{n}^{\xi}(s)))^{*})
  \notag \\&
 -\int_{H\backslash\left\{0\right\}}\left[\Psi(s, X_{n}^{\xi}(s)+R_{n}f(v,X_{n}^{\xi}(s)))-\Psi(s,X_{n}^{\xi}(s))-\left\langle \Psi_{x}(s,X_{n}^{\xi}(s)),R_{n}f(v,X_{n}^{\xi}(s))\right\rangle_{H}\right]\beta(dv)|.
 \notag\\&
 =|\lim_{n\rightarrow\infty}\left\langle \Psi_{x}(s,X_{n}^{\xi}(s)),AX_{n}^{\xi}(s)\right\rangle
       -\lim_{n\rightarrow\infty}\left\langle \Psi_{x}(s,X_{n}^{\xi}(s)),AX_{n}^{\xi}(s)\right\rangle|=0  \quad P-a.s. .\notag
\end{align}
Similarly as in Theorem \ref{imldtm4} the r.h.s. of \eqref{limL} can be divided in terms which modulus converges to zero a.s.. This is proven using  Theorem \ref{imldtm3}, the fact that  
 $\left\|R_n\right\|$ is uniformly bounded, as well as $(R_n-I)x\rightarrow 0$, that 
 \eqref{sup} holds,  $E[\sup_{0\leq s\leq T}\left\|X^\xi(s)\right\|_{H}^{2}]<\infty$ and conditions \eqref{bd-h-1}, \eqref{bd-h-2}, \eqref{bd-h-3} hold.
\end{proof}
\end{corollary}

Here we will present an example of Theorem \ref{imldtm4}.\\

\noindent\textbf{Example 1:-}  \\
 Let $\Lambda :H\rightarrow\mathbb{R}$. Assume that the assumptions of Theorem \ref{imldtm4} hold. Now if we apply 
 the It\^o formula for mild solutions, of \eqref{imldeq1}, for the function $e^{ct}\Lambda(x)$ where $c>0$, $t\geq0$ then it will be   
\begin{align*}
 &\lim_{n\rightarrow\infty}\int_{0}^{t}e^{cs}\left\langle\Lambda'(X_{n}^{\xi}(s)),AX_{n}^{\xi}(s)\right\rangle_{H}ds&\\
 &=e^{ct}\Lambda(X^{\xi}(t))-\Lambda(\xi)-\int_{0}^{t}ce^{cs}(\Lambda(X^{\xi}(s)))ds&\\
 &-\int_{0}^{t}e^{cs}\left\langle\Lambda'(X^{\xi}(s)),F(X^{\xi}(s))\right\rangle_{H}ds
       -\int_{0}^{t}e^{cs}\frac{1}{2} tr(\Lambda''(X^{\xi}(s))(B(X^{\xi}(s)))Q(B(X^{\xi}(s)))^{*})ds&\\
 &-\int_{0}^{t}\int_{H\backslash\left\{0\right\}}e^{cs}\left[\Lambda(X^{\xi}(s)+f(v,X^{\xi}(s)))
      -\Lambda(X^{\xi}(s))-\left\langle \Lambda'(X^{\xi}(s)),f(v,X^{\xi}(s))\right\rangle_{H}\right]\beta(dv)ds&\\
 &-\int_{0}^{t}e^{cs}\left\langle \Lambda'(X^{\xi}(s)),B(X^{\xi}(s))dW_{s}\right\rangle_{H}&\\
 &-\int_{0}^{t}\int_{H\backslash\left\{0\right\}}e^{cs}\left[\Lambda(X^{\xi}(s)+f(v,X^{\xi}(s)))-\Lambda(X^{\xi}(s))\right]q(dv,ds)&
\end{align*}
 $P$-a.s..  
 
 \begin{remark} \label{Rem5}  By Corollary \ref{Rem3} we can take in Example 1   e.g.  $\Lambda \in C_b^{2}(\mathbb{R}\times H)$, but it can also be a function in $C^{2}(\mathbb{R}\times H)$ for which \eqref{bd-h-1},   \eqref{bd-h-2} and   \eqref{bd-h-3} are satisfied. Suppose  e.g. 
  $\int_{H\setminus \{0\}} (\|f(v,x)\|^2+\|f(v,x)\|^4) \beta(dv) <\infty $, for all $x\in H$  and then choose e.g. $\Lambda(x)=\|x\|^2$ (see the examples in \cite{9}, page 184).
 \end{remark}

 
   \subsection{Applications of the It\^o formula obtained by Yosida approximations: Stability properties of mild solutions}\label{Applications stability}
   In this section we will present  how the It\^o formula obtained by Yosida approximation provides a good method to prove that the mild solution of \eqref{imldeq1} is  "exponential stable  in the   mean square sense" ("exponential stable in the m.s.s") or "exponential ultimate bounded in the  mean square sense"  ("exponential ultimate bounded in the m.s.s.") under suitable conditions. In the whole section we assume that conditions (A1), (A2) and (A3) hold. Then due to Theorem \ref{imldtm1}
      there exists a unique    mild solution $\{X^x(t), t\geq 0\}$  of \eqref{imldeq1} with initial condition $X^x(0)=x$$\in H$.\\

   \begin{definition}\label{m.s.s-ex} We say that that the solution of \eqref{imldeq1} is exponentially stable in the 
   mean square sense (m.s.s.) if there exist $c,\beta>0$ such that for all $t\geq 0$ 
   and $x\in H$,
   
   \begin{align}\label{m.s.s.}
   E\left\|X^x(t)\right\|^2_H\leq ce^{-\beta t}\left\|x\right\|^2_H; .
   \end{align}
   \end{definition}

Here as usual $C^2(H)$ denotes the space of continuous functions  $\Psi:H\rightarrow\mathbb{R}$, with continuous
 first and second  Fr\'echet  derivative $\Psi'(x)$ and $\Psi''(x)$, $x\in H$.

 \begin{definition}\label{lypu-fn} Let  $\mathcal{L}$ be defined as  in (\ref{defnL}).
   A function $\Psi:H\rightarrow\mathbb{R} \in C^2(H)$ is called a Lyapunov function for $\mathcal{L}$ if it satisfies  the following conditions:\\
   
   (I) There exist finite constants $c_1$, $c_2>0$ such that; for all $x\in H$
   
   \begin{align*}
   c_1\left\|x\right\|_H^2\leq\Psi(x)\leq c_2\left\|x\right\|_H^2 
   \end{align*}
   
   
   (II) There exists a constants $c_3>0$ such that
   \begin{align*}
   \mathcal{L}\Psi(x)\leq-c_3\Psi(x)\ \ \ \textit{for all}\ \ \ x\in\mathcal{D}(A)
   \end{align*}
  .

   \end{definition}

 \begin{theorem}\label{m.s.s-th1} Assume there exists a  function $\Psi \in   C^2(H)$ which is a Lyapunov function  for $\mathcal{L}$, and  such that \eqref{bd-h-1},   \eqref{bd-h-2} and   \eqref{bd-h-3} are satisfied, then 
 the mild solution of \eqref{imldeq1} is exponentially stable in the m.s.s.

 \end{theorem}
 
 If we assume that the first two  Frech\'et derivatives of the Lyapunov function $\Psi$ are uniformly bounded, i.e. $\Psi \in   C_b^2(H)$,  then due to Remark  \ref{Rem1} it follows directly that the following Corollary holds:
  \begin{corollary}
   \label{corstablemssu}
    
    If there exists a  function $\Psi \in   C_b^2(H)$ which is a Lyapunov function for $\mathcal{L}$, then the mild solution of \eqref{imldeq1}   is exponentially stable in the m.s.s. 
   \end{corollary}

  A detailed proof of Corollary \ref{corstablemssu} is in Theorem 6.4 of \cite{1}, Theorem 4.2 of \cite{3} (see also Section 7.1 in \cite{10}). Here we present a simplified proof of the more general statement in Theorem \ref{m.s.s-th1} by using the  limiting 
  argument of Theorem \ref{imldtm4} and Corollary \ref{itocorollary}.\\
 \begin{proof}
  We first prove that inequality \eqref{m.s.s.} holds for all $x\in \mathcal{D}(A)$.\\
 By applying the It\^o formula to the Yosida approximation $X_n^x$ of $X^x$  and by  taking expectations on both sides we obtain
 
 \begin{align}\label{recorec1}
 e^{c_3t}E\Psi(X_n^x(t))-\Psi(X_n^x(0))=E\int_0^te^{c_3s}\left(c_3\Psi(X_n^x(s))+\mathcal{L}_n\Psi(X_n^x(s))\right)ds.\quad \forall \, x\in \mathcal{D}(A). 
 \end {align}
 
 From condition (II),
 
 \begin{align*} 
 c_3\Psi(X_n^x(s))+\mathcal{L}_n\Psi(X_n^x(s))\leq-\mathcal{L}\Psi(X_n^x(s))+\mathcal{L}_n\Psi(X_n^x(s))\quad \forall \, x\in \mathcal{D}(A). 
 \end {align*}
 \begin{align}\label{recorec2}
 \Rightarrow e^{c_3t}E\Psi(X_n^x(t))-\Psi(X_n^x(0))\leq E\int_0^te^{c_3s}\left(-\mathcal{L}\Psi(X_n^x(s))+\mathcal{L}_n\Psi(X_n^x(s))\right)ds.\quad \forall \, x\in \mathcal{D}(A). 
 \end {align}
 
 Using Corollary \ref{itocorollary} and \eqref{sup} we obtain  
 
 \begin{align}\label{recorec4}
 e^{c_3t}E\Psi(X^x(t))\leq\Psi(x)
 \end{align}
 \begin{align*}
 \Rightarrow c_1E\left\|X^x(t)\right\|^2_H\leq E\Psi(X^x(t))\leq e^{-c_3t}\Psi(x)\leq c_2e^{-c_3t}\left\|x\right\|^2_H \ \ \ \ 
 [\text{from condition (ii)}] 
 \end{align*}
 \begin{align}\label{recorec5}
 \Rightarrow E\left\|X^x(t)\right\|^2_H\leq \frac{c_2}{c_1}e^{-c_3t}\left\|x\right\|^2_H \forall \, x\in \mathcal{D}(A). 
 \end{align}
 Choosing $c=\frac{c_2}{c_1}$ and $\beta=c_3$, we can conclude that \eqref{m.s.s.} holds for all $x\in \mathcal{D}(A)$.
 As however $\mathcal{D}(A)$ is dense in $H$ and, according to Corollary 5.3.2 in \cite{10}, for each $t>0$ there exists a constant $C_t$ such that 
 \begin{equation}\label{ineqindata}
 \mathbb{E}\left\|X^x(t)-X^y(t)\right\|^2_H\leq C_t \left\|x-y\right\|^2_H,
 \end{equation}
 we can conclude that, the mild solution $X^x(t)$ is exponentially stable in 
   the mean square sense (m.s.s.) 
 
 ..
 \end{proof}

\begin{definition}\label{ulti-bnDfn}
 We say that the mild solution of (\eqref{imldeq1}) is exponentially ultimate bounded in the mean square sense (m.s.s.)  
 if there exist positive constants $c$, $\beta$, $M$ such that 
\begin{align}\label{u.b.m.s.s.}
E\left\|X^x(t)\right\|^2_H\leq ce^{-\beta t}\left\|x\right\|^2_H+M;\ \ \ \text{for all}\ x\in H \text{and} t> 0.
\end{align}
 \end{definition}
 
 \begin{theorem}\label{ulti-bn-th1}
 Let us assume there exists a function $\Psi\in C^2(H)$ satisfying the following conditions: \\
 
 (i) $c_1\left\|x\right\|_H^2-k_1\leq\Psi(x)\leq c_2\left\|x\right\|_H^2-k_2$; for all $x\in H$\\ 
 
 
 (ii) $\mathcal{L}\Psi\leq-c_3\Psi(x)+k_3$; for $x\in\mathcal{D}(A)$, \\
 
 \noindent where $c_1$, $c_2$, $c_3$, $k_1$, $k_2$ and $k_3$ are finite, positive constants.\\
  Moreover, let us assume 
  \eqref{bd-h-1},   \eqref{bd-h-2} and   \eqref{bd-h-3} are satisfied.\\ 
  Then the mild solution $\{X^x(t),t\geq 0\}$ of \eqref{imldeq1} 
  is exponentially ultimate bounded in the m.s.s.

 \end{theorem}
  By similar arguments  as in the proof of Corollary \ref{corstablemssu} the following Corollary follows directly
    \begin{corollary}\label{corubddmssu}
        If there exists a  function $\Psi \in   C_b^2(H)$ such that conditions (i) and (ii) hold in Theorem \ref{ulti-bn-th1}, then the mild solution $\{X^x(t),t\geq 0\}$ of \eqref{imldeq1} is exponentially ultimate bounded in the m.s.s..
      \end{corollary}

  A detailed proof of Corollary \ref{corubddmssu} is in Theorem 7.1 of \cite{1}, Theorem 5.2 of \cite{3} (see also Section 7.2 in \cite{10}). Here we present a simplified proof of the more general statement in Theorem \ref{ulti-bn-th1} by using the  limiting 
  argument of Theorem \ref{imldtm4} and Corollary \ref{itocorollary}.\\
 \begin{proof} First let us assume $x\in \mathcal{D}(A)$.By applying the It\^o formula to the Yosida approximation $X_n^x$ of $X^x$  and by  taking expectations on both sides we obtain 
 
\begin{align}\label{recorec100}
e^{c_3t}E\Psi(X_n^x(t))-\Psi(X_n^x(0))=E\int_0^te^{c_3s}\left(c_3\Psi(X_n^x(s))+\mathcal{L}_n\Psi(X_n^x(s))\right)ds.
\end {align}

Now from condition (ii) it follows 

\begin{align*} 
c_3\Psi(X_n^x(s))+\mathcal{L}_n\Psi(X_n^x(s))\leq-\mathcal{L}\Psi(X_n^x(s))+k_3+\mathcal{L}_n\Psi(X_n^x(s)).
\end {align*}
 \begin{align}\label{recorec101}
 \Rightarrow e^{c_3t}E\Psi(X_n^x(t))-\Psi(X_n^x(0))&\leq E\int_0^te^{c_3s}\left(-\mathcal{L}\Psi(X_n^x(s))+k_3+\mathcal{L}_n\Psi(X_n^x(s))\right)ds&\\
&=E\int_0^te^{c_3s}\left(-\mathcal{L}\Psi(X_n^x(s))+\mathcal{L}_n\Psi(X_n^x(s))\right)ds+\int_0^te^{c_3s}k_3ds.&\nonumber
\end {align}

Due to Corollary \ref{itocorollary} we obtain

\begin{align*}
E\int_0^te^{c_3s}\left(-\mathcal{L}\Psi(X_n^x(s))+\mathcal{L}_n\Psi(X_n^x(s))\right)ds\rightarrow 0.
\end{align*}



Therefore as $n\rightarrow\infty$,  using the continuity of $\Psi$ and Lebesgue dominated convergence Theorem, we get

\begin{align}\label{recorec102}
e^{c_3t}E\Psi(X^x(t))&\leq\Psi(x)+\int_0^te^{c_3s}k_3ds&\\
&=\Psi(x)+\frac{k_3}{c_3}(e^{c_3t}-1).& \nonumber
\end{align}
\begin{align}\label{recorec103}
\Rightarrow E\Psi(X^x(t))\leq e^{-c_3t}\Psi(x)+\frac{k_3}{c_3}(1-e^{-c_3t}).
\end{align}

Now from condition (i) we get that for all $x\in H$

\begin{align}\label{pg4}
 c_1E\left\|X^x(t)\right\|_H^2-k_1\leq E\Psi(X^x(t))&\leq e^{-c_3t}\left(c_2\left\|x\right\|_H^2-k_2\right)+\frac{k_3}{c_3}(1-e^{-c_3t})&\\
&\leq c_2e^{-c_3t}\left\|x\right\|_H^2+\frac{k_3}{c_3}(1-e^{-c_3t}).\nonumber&
\end{align}
\begin{align*}
\Rightarrow c_1E\left\|X^x(t)\right\|_H^2\leq c_2e^{-c_3t}\left\|x\right\|_H^2+\frac{k_3}{c_3}(1-e^{-c_3t})+k_1,
\end{align*}
\begin{align*}
\Rightarrow E\left\|X^x(t)\right\|_H^2\leq\frac{c_2}{c_1}e^{-c_3t}\left\|x\right\|_H^2+\frac{1}{c_1}\left(k_1+\frac{k_3}{c_3}\right).
\end{align*}
By using the density of $\mathcal{D}(A)$ in $H$ and inequality  \eqref{ineqindata}  we conclude that by 
 choosing $c=\frac{c_2}{c_1}$, $\beta=c_3$ and $M=\frac{1}{c_1}\left(k_1+\frac{k_3}{c_3}\right)$ inequality \eqref{u.b.m.s.s.} holds for all $x\in H$,  and hence the mild solution $\{X^x(t),t\geq 0\}$ of \eqref{imldeq1} is exponentially ultimate bounded in the m.s.s..

 \end{proof}

 \subsection {An "Ichikawa-type" It\^o  formula for mild solutions}  \label{ichikawa}

  In this subsection  we prove that the Ichikawa type It\^o formula for the mild solutions  $\left\{X^{\xi}(t),t\geq0\right\}$ obtained by Ichikawa for SPDEs driven by 
 Gaussian noise in \cite{8}, can also be generalized to the case of the  SPDE defined in (\ref{imldeq1}).
The advantage of introducing this "Ichikawa -type" It\^o  formula  for mild solutions of (\ref{imldeq1}) is that in  this type of  It\^o  formula (\eqref{ichiIto1}) below) there is no Yosida approximant appearing, 
unlike than in the It\^o  formula \eqref{ito-lim-forml} obtained in Theorem \ref{imldtm4}. The disadvantage is however that compared to Theorem \ref{imldtm4} we have to restrict the set of functions $\Psi: [0,T]\times H\rightarrow\mathbb{R}$ for which we can apply the It\^o formula  To this aim we introduce the following notation:

\par 
  Let $\tilde{C}^{1,2}([0,T]\times H)$ (resp. $\tilde{C}_b^{1,2}([0,T]\times H)$) be the class of functions $\Psi\in C^{1,2}([0,T]\times H)$ (resp. $\Psi\in C_b^{1,2}([0,T]\times H)$) with the properties: \\
 (I1) The function $\mathcal{L}\Psi(s,x)$ can be extended to a continuous function $\overline{\mathcal{L}\Psi}(s,x)$ on $[0,T]\times H$ for $x\in H$.  \\
 (I2) $\|\overline{\mathcal{L}\Psi}(s,x)\|\leq k(1+\|x\|^{2})$, for $x\in H$, 
  $s\in[0,T]$ and for some $k>0$. \\

  \begin{theorem}\label{ichiIto}  Let $\{S(t),t\geq0\}$ be a pseudo-contraction semigroup generated by $A$ satisfying assumption (A1). Suppose assumption (A2) holds. 
  Let the coefficients $F$, $B$, $f$ satisfy assumption (A3).    Let  $\{X^\xi(t),t\in [0,T]\}$ be the  mild solution of (\ref{imldeq1}) in $[0,T]$ with initial condition $\xi \in \mathcal{D}(A)$.
 For any $\Psi$ $\in \tilde{C}^{1,2}([0,T]\times H)$   such that \eqref{bd-h-1}, \eqref{bd-h-2},  \eqref{bd-h-3} are satisfied, the following Ichikawa type It\^o Formula for mild solutions hold $P$-a.s. for all $t\in[0,T]$ 
 
 \begin{align}\label{ichiIto1}
 &\Psi(t,X^{\xi}(t))-\Psi(0,\xi)=\int_{0}^{t}(\Psi_{s}(s,X^{\xi}(s))+\overline{\mathcal{L}\Psi}(s,X^{\xi}(s)))ds&\\
 &+\int_{0}^{t}\left\langle \Psi_{x}(s,X^{\xi}(s)),B(X^{\xi}(s))dW_{s}\right\rangle_{H}&\nonumber\\
 &+\int_{0}^{t}\int_{H\backslash\left\{0\right\}}\left[\Psi(s,X^{\xi}(s)+f(v,X^{\xi}(s)))-\Psi(s,X^{\xi}(s))\right]q(dv,ds).&\nonumber
 \end{align}

 \begin{proof}
  Here we assumed that $\Psi\in\tilde{C}^{1,2}([0,T]\times H)$. So, it satisfies the conditions (I1) and (I2). 
  From (I1) the continuous extension 
  of $\mathcal{L}\Psi(s,x)$ exists in $[0,T]\times H$, which is $\overline{\mathcal{L}\Psi}(s,x)$.    It follows that if we denote again with    $\left\{X_{n}(t)^{\xi},t\geq 0\right\}$ the Yosida approximation of   $\left\{X^{\xi}(t),t\geq 0\right\}$, then 
  
  \begin{align*}
&  \overline{\mathcal{L}\Psi}(s,X^{\xi}(s))
  =\lim_{n\rightarrow\infty}\mathcal{L}\Psi(s,X_{n}^{\xi}(s)) \\
  &=\left\langle \Psi_{x}(s,X^{\xi}(s)),F(X^{\xi}(s))\right\rangle_{H}
         +\lim_{n\rightarrow\infty}\left\langle\Psi_{x}(s,X_{n}^{\xi}(s)),AX_{n}^{\xi}(s)\right\rangle_{H}&  \\
  &+ \frac{1}{2} tr(\Psi_{xx}(s,X^{\xi}(s))(B(X^{\xi}(s)))Q(B(X^{\xi}(s)))^{*})& \\ 
  &+\int_{H\backslash\left\{0\right\}}\left[\Psi(s,X^{\xi}(s)+f(v,X^{\xi}(s)))-\Psi(s,X^{\xi}(s))
         -\left\langle \Psi_{x}(s,X^{\xi}(s)),f(v,X^{\xi}(s))\right\rangle_{H}\right]\beta(dv)&
  \end{align*} 
   where $X_{n}^{\xi}(s)\in\mathcal{D}(A)$. We know here that the limit  $\lim_{n\rightarrow\infty}\left\langle\Psi_{x}(s,X_{n}^{\xi}(s)),AX_{n}^{\xi}(s)\right\rangle_{H}$ is well defined, as we assume (I1) and all other terms converge, as shown in the proof of Theorem \ref{imldtm4}. 
  By (I2) $\overline{\mathcal{L}\Psi}(s,X^{\xi}(s))$ is bounded by integrable function, so by applying Lebesgue dominated convergence theorem, we  conclude that 
  $\int_0^t\mathcal{L}\Psi(s,X_{n}^{\xi}(s))ds\rightarrow\int_0^t\overline{\mathcal{L}\Psi}(s,X^{\xi}(s))ds$, as $n\to \infty$. Hence we  conclude that the Ichikawa type It\^o Formula of (\ref{ichiIto1}) holds.
 \end{proof}
 \end{theorem}
    
 The  proof of following Corollary is straight by  using   Remark \ref{Rem2}.
 \begin{corollary} \label{ichiIto2}
  Let $\{S(t),t\geq 0\}$ be a pseudo-contraction semigroup generated by $A$ satisfying assumption (A1). Suppose assumption (A2) holds. 
    Let the coefficients $F$, $B$, $f$ satisfy assumption (A3).   Let  $\{X^\xi(t),t\in [0,T]\}$ be the  mild solution of (\ref{imldeq1}) in $[0,T]$ with initial condition  $\xi \in \mathcal{D}(A)$.
   Then for all $\Psi$ $\in \tilde{C}_b^{1,2}([0,T]\times H)$  the Ichikawa -type It\^o formula \eqref{ichiIto1} holds.
 \end{corollary}

 \begin{remark}\label{Remark 2} Whenever $\Psi\in\tilde{C}^{1,2}([0,T]\times H)$ ($\Psi$ $\in \tilde{C}_b^{1,2}([0,T]\times H)$) i.e. 
 $\overline{\mathcal{L}\Psi}(s,x)$ exists and the conditions of Theorem \ref{ichiIto} (Theorem \ref{ichiIto2}) are satisfied, then we can interchange the limit with the 
 integral in the It\^o formula for mild solutions of Theorem \ref{imldtm4}. Then the It\^o formula for mild solutions of Theorem \ref{imldtm4} 
 (eq. \eqref{ito-lim-forml}) can be rewritten as follows-
 \begin{align}\label{ito-lim-forml-rmrk}
 &\int_{0}^{t}\lim_{n\rightarrow\infty}\left\langle\Psi_{x}(s,X_{n}^{\xi}(s)),AX_{n}^{\xi}(s)\right\rangle_{H}ds&\\
 &=\Psi(t,X^{\xi}(t))-\Psi(0,\xi)\nonumber-\int_{0}^{t}(\Psi_{s}(s,X^{\xi}(s)))ds
 -\int_{0}^{t}\left\langle\Psi_{x}(s,X^{\xi}(s)),F(X^{\xi}(s))\right\rangle_{H}ds&\nonumber\\
 &-\int_{0}^{t}\frac{1}{2} tr(\Psi_{xx}(s,X^{\xi}(s))(B(X^{\xi}(s)))Q(B(X^{\xi}(s)))^{*})ds&\nonumber\\
 &-\int_{0}^{t}\int_{H\backslash\left\{0\right\}}\left[\Psi(s,X^{\xi}(s)+f(v,X^{\xi}(s)))-\Psi(s,X^{\xi}(s))
  -\left\langle \Psi_{x}(s,X^{\xi}(s)),f(v,X^{\xi}(s))\right\rangle_{H}\right]\beta(dv)ds&\nonumber\\
 &-\int_{0}^{t}\left\langle \Psi_{x}(s,X^{\xi}(s)),B(X^{\xi}(s))dW_{s}\right\rangle_{H}&\nonumber\\
 &-\int_{0}^{t}\int_{H\backslash\left\{0\right\}}\left[\Psi(s,X^{\xi}(s)+f(v,X^{\xi}(s)))-\Psi(s,X^{\xi}(s))\right]q(dv,ds).&\nonumber
\end{align}
\end{remark}

 \vspace{0.5cm}
 Here we show through two examples how the It\^o formula \eqref{ito-lim-forml} in Theorem \ref{imldtm4}, obtained through Yosida approximation,  can be related  to the "Ichikawa type" -It\^o formula \eqref{ichiIto1} in Theorem \ref{ichiIto}.\\

  \noindent\textbf{Example 2:-}  \\
  Let $A$ be a symmetric linear operator and the assumptions of Theorem \ref{imldtm4} hold. Assume that, for a fixed $l>0$, $l\in(0,\infty)$, 
  $\Psi(s,x)=e^{lA}\varLambda(s,x)$ and $\varLambda_x(s,x)\in\mathcal{D}(A)$.  
 Also assume that 
  $\Psi$$\in C^{1,2}(\mathbb{R_+}\times H)$ and  satisfies the following property,  \\
  (i) There exists constants $c_{1}>0$   s.t.
$\|\Psi_{x}(s,x)\|\leq c_{2}\left\|x\right\|_{H};$
  for all $x\in H$.  \\
  
  Moreover assume  \eqref{bd-h-1}, \eqref{bd-h-2},  \eqref{bd-h-3} are satisfied.\\
  
  \begin{align}\label{ex4}
 \left\langle\Psi_{x}(s,X_{n}^{\xi}(s)),AX_{n}^{\xi}(s)\right\rangle_{H}\leq\mathcal{C}\left\|X_{n}^{\xi}(s)\right\|_{H}^{2}.
 \end{align}


  Here,
 
 \begin{align*}
 &\left\langle\Psi_{x}(s,X_{n}^{\xi}(s)),AX_{n}^{\xi}(s)\right\rangle_{H}&\\
 &=\left\langle e^{lA}\varLambda_{x}(s,X_{n}^{\xi}(s)),AX_{n}^{\xi}(s)\right\rangle_{H}&\\
 &=\left\langle A^{*}e^{lA}\varLambda_{x}(s,X_{n}^{\xi}(s)),X_{n}^{\xi}(s)\right\rangle_{H}&\\
 &\leq\mathcal{C}\left\|X_{n}^{\xi}(s)\right\|_{H}^{2}.&
 \end{align*}
 We can write from second to third equality, due to the fact that,  
  $e^{lA}\varLambda_x(s,X_{n}^{\xi}(s))\in\mathcal{D}(A)\subset\mathcal{D}(A^*)$ as we assumed $\varLambda_x(s,x)\in\mathcal{D}(A)$ 
  and $A$ to be a symmetric linear operator.

   By Theorem \ref{imldtm3} and Generalized Lebesgue Dominated Convergence theorem (Theorem 3.4 of \cite{1}) we  conclude that,
    \begin{align*}
  \int_{0}^{t}\left\langle e^{lA}\varLambda_{x}(s,X_{n}^{\xi}(s)),AX_{n}^{\xi}(s)\right\rangle_{H}ds\rightarrow\int_{0}^{t}\left\langle A^{*}e^{lA}\varLambda_{x}(s,X^{\xi}(s)),X^{\xi}(s)\right\rangle_{H}ds
  \end{align*}
  
  $P$-a.s.. Therefore we can write the It\^o formula for mild solutions of Theorem \ref{imldtm4} for the function 
  $\Psi(s,x)=e^{lA}\varLambda(s,x)$, for $x\in H$, as follows
  
  \begin{align}\label{exampleIt}
   &e^{lA}\varLambda(t,X^{\xi}(t))-e^{lA}\varLambda(0,\xi)=\int_{0}^{t} e^{lA}\varLambda_{s}(s,X^{\xi}(s))ds&\\
   &+\int_{0}^{t}\left\langle e^{lA}\varLambda_{x}(s,X^{\xi}(s)),F(X^{\xi}(s))\right\rangle_{H}ds
       +\int_{0}^{t}\left\langle A^{*}e^{lA}\varLambda_{x}(s,X^{\xi}(s)),X^{\xi}(s)\right\rangle_{H}ds&\nonumber\\
  &+\int_{0}^{t}\frac{1}{2} tr(e^{lA}\varLambda_{xx}(s,X^{\xi}(s))(B(X^{\xi}(s)))Q(B(X^{\xi}(s)))^{*})ds&\nonumber\\
  &+\int_{0}^{t}\int_{H\backslash\left\{0\right\}}e^{lA}\big[\varLambda(s,X^{\xi}(s)+f(v,X^{\xi}(s)))-\varLambda(s,X^{\xi}(s))
          -\left\langle\varLambda_{x}(s,X^{\xi}(s)),f(v,X^{\xi}(s))\right\rangle_{H}\big]\beta(dv)ds&\nonumber\\
  &+\int_{0}^{t}\left\langle e^{lA}\varLambda_{x}(s,X^{\xi}(s)),B(X^{\xi}(s))dW_{s}\right\rangle_{H}&\nonumber\\
  &+\int_{0}^{t}\int_{H\backslash\left\{0\right\}}e^{lA}\left[\varLambda(s,X^{\xi}(s)+f(v,X^{\xi}(s)))-\varLambda(s,X^{\xi}(s))\right]q(dv,ds)&\nonumber
  \end{align}
  $P$-a.s. and by Theorem \ref{ichiIto}, we obtain 
  \begin{align*}
  &\overline{\mathcal{L}\Psi}(s,X^{\xi}(s))=\left\langle e^{lA}\varLambda_{x}(s,X^{\xi}(s)),F(X^{\xi}(s))\right\rangle_{H}
        +\left\langle A^{*}e^{lA}\varLambda_{x}(s,X^{\xi}(s)),X^{\xi}(s)\right\rangle_{H}&\\
  &+\frac{1}{2} tr(e^{lA}\varLambda_{xx}(s,X^{\xi}(s))(B(X^{\xi}(s)))Q(B(X^{\xi}(s)))^{*})&\\
  &+\int_{H\backslash\left\{0\right\}}e^{lA}\left[\varLambda(s,X^{\xi}(s)+f(v,X^{\xi}(s)))-\varLambda(s,X^{\xi}(s))
        -\left\langle\varLambda_{x}(s,X^{\xi}(s)),f(v,X^{\xi}(s))\right\rangle_{H}\right]\beta(dv).&
  \end{align*}
  
  \vspace{0.5cm}
   
   \noindent\textbf{Example 3:-}  \\
   Let $A$ be a symmetric linear operator. 
   Now, we consider $\Psi(s,x)=e^{(t-s)A}\varGamma(x)$ for $x\in H$. Assume that  $\Psi$$\in C^{1,2}(\mathbb{R_+}\times H)$ satisfies the condition (i) of Example 2. Also assume that 
   $\varGamma, \varGamma_x\in\mathcal{D}(A)$.  Moreover assume  \eqref{bd-h-1}, \eqref{bd-h-2},  \eqref{bd-h-3} are satisfied.\\

   \begin{align*}
   \Psi_{s}(s,X_{n}^{\xi}(s))=\left(-Ae^{(t-s)A}\varGamma(X_{n}^{\xi}(s))\right)\rightarrow\left(-Ae^{(t-s)A}\varGamma(X^{\xi}(s))\right)
   \end{align*}
   
   $P$-a.s.. The convergence of all other terms are similar as in Example 2. So, we can write the It\^o formula of Theorem \ref{imldtm4} as,

   \begin{align}\label{rocex}
   &\varGamma(X^{\xi}(t))-e^{tA}\varGamma(\xi)=\int_{0}^{t}\left(-Ae^{(t-s)A}\varGamma(X^{\xi}(s))\right)ds&\\
   &+\int_{0}^{t}\left\langle e^{(t-s)A}\varGamma_{x}(X^{\xi}(s)),F(X^{\xi}(s))\right\rangle_{H}ds
           +\int_{0}^{t}\left\langle A^{*}e^{(t-s)A}\varGamma_{x}(X^{\xi}(s)),X^{\xi}(s)\right\rangle_{H}ds&\nonumber\\
   &+\int_{0}^{t}\frac{1}{2} tr(e^{(t-s)A}\varGamma_{xx}(X^{\xi}(s))(B(X^{\xi}(s)))Q(B(X^{\xi}(s)))^{*})ds&\nonumber\\
   &+\int_{0}^{t}\int_{H\backslash\left\{0\right\}}[e^{(t-s)A}\varGamma(X^{\xi}(s)+f(v,X^{\xi}(s)))-e^{(t-s)A}\varGamma(X^{\xi}(s))&\nonumber\\
   &\ \ \ \ \ \ \ \ \ \ \ \ \ \ \ \ \ \ -\left\langle e^{(t-s)A}\varGamma_{x}(X^{\xi}(s)),f(v,X^{\xi}(s))\right\rangle_{H}]\beta(dv)ds&\nonumber\\
   &+\int_{0}^{t}\left\langle e^{(t-s)A}\varGamma_{x}(X^{\xi}(s)),B(X^{\xi}(s))dW_{s}\right\rangle_{H}&\nonumber\\
   &+\int_{0}^{t}\int_{H\backslash\left\{0\right\}}\left[e^{(t-s)A}\varGamma(X^{\xi}(s)+f(v,X^{\xi}(s)))
          -e^{(t-s)A}\varGamma(X^{\xi}(s))\right]q(dv,ds)&\nonumber
   \end{align}

   $P$-a.s. and by Theorem \ref{ichiIto}, we obtain 
   
   \begin{align*}
   &\overline{\mathcal{L}\Psi}(s,X^{\xi}(s))=\left\langle e^{(t-s)A}\varGamma_{x}(X^{\xi}(s)),F(X^{\xi}(s))\right\rangle_{H}
         +\left\langle A^{*}e^{(t-s)A}\varGamma_{x}(X^{\xi}(s)),X^{\xi}(s)\right\rangle_{H}&\\
   &+\frac{1}{2} tr(e^{(t-s)A}\varGamma_{xx}(X^{\xi}(s))(B(X^{\xi}(s)))Q(B(X^{\xi}(s)))^{*})&\\
   &+\int_{H\backslash\left\{0\right\}}[e^{(t-s)A}\varGamma(X^{\xi}(s)+f(v,X^{\xi}(s)))-e^{(t-s)A}\varGamma(X^{\xi}(s))&\\
   &\ \ \ \ \ \ \ \ \ \ \ \ \ \ \ \ \ \ -\left\langle e^{(t-s)A}\varGamma_{x}(X^{\xi}(s)),f(v,X^{\xi}(s))\right\rangle_{H}]\beta(dv).&
   \end{align*}
   
  \vspace{0.5cm}
   
  \section{An It\^o formula type inequality}\label{relt-ichi-sem-gen}
 
  We assume in the whole section that  
  $\{S(t),t\geq0\}$ is a pseudo-contraction semigroup generated by $A$ satisfying assumption (A1), that (A2) holds. 
 and  that the coefficients $F$, $B$, $f$ satisfy assumption (A3).
  Let  $X^{x}(t):==X(t,0;x)\in H$ denote the  the  mild solution of \eqref{imldeq4}
 with deterministic initial condition $x\in H$. $X^{x}(t)$ is a homogeneous Markov process  (see, section 3.4 of \cite{1} and section 6 of \cite{4}), section 5.4 \cite{10}. \\
Let  $\Psi\in C_b^{2}(H)$ and  $x\in\mathcal{D}(A)$, then the   the operator $\mathcal{L}$ defined in \eqref{defnL} is well defined in all its terms and acts on $\Psi$ as follows:  
 . 
 \begin{align*}
 \mathcal{L}\Psi(x):=&\left\langle \Psi_{x}(x),Ax+F(x)\right\rangle_{H}+\frac{1}{2} tr(\Psi_{xx}(x)(B(x))Q(B(x))^{*})&\\
 &+\int_{H\backslash\left\{0\right\}}\left[\Psi(x+f(v,x))-\Psi(x)-\left\langle \Psi_{x}(x),f(v,x)\right\rangle_{H}\right]\beta(dv).&
 \end{align*}

Let $P_{t}$ be the semigroup acting on  $\Psi\in C_{b}^{2}(H)$, such that 
  \begin{align*}
  [P_{t}\Psi](x)=E[\Psi(X^{x}(t))]\quad for x\in H
  \end{align*}

  \begin{theorem}\label{ichigenerator}
  Assume that the solution $X^{x}(t)$ of SPDE \eqref{imldeq1} with initial condition $x\in \mathcal{D}(A)$ satisfies  $X^x(t)\in\mathcal{D}(A)$ for all $t \geq 0$. Then for  $\Psi\in C_{b}^{2}(H)$ 
  \begin{align}\label{relation1}
  [P_{t}\Psi](x)-\Psi(x)=\int_{0}^{t}[P_{s}\mathcal{L}\Psi](x)ds
  \end{align}
  and 
  \begin{align}\label{relation2}
  \lim_{t\downarrow0}\frac{[P_{t}\Psi](x)-\Psi(x)}{t}=[\mathcal{L}\Psi](x).
  \end{align}
  \begin{proof} 
   We refer sections 4.1, 4.2 of \cite{kurtz} or sections 3.2, 3.3 of \cite{6} for related theory.  
     \begin{align}\label{imldeq6t}
    \Psi(X^{x}(t))-\Psi(x)=\int_{0}^{t}\mathcal{L}\Psi(X^{x}(s))ds
    \end{align}
    \begin{align*}
    +\int_{0}^{t}\left\langle \Psi_{x}(X^{x}(s)),B(X^{x}(s))dW_{s}\right\rangle_{H}
    \end{align*}
    \begin{align*}
    +\int_{0}^{t}\int_{H\backslash\left\{0\right\}}\left[\Psi(X^{x}(s)+f(v,X^{x}(s)))-\Psi(X^{x}(s))\right]q(dv,ds).
    \end{align*}
    Now take the expectation on both sides of (\ref{imldeq6t}). Second and third term of R.H.S. of (\ref{imldeq6t}) will be zero, because of martingale and 
    we get,
    \begin{align*}
    E[\Psi(X^{x}(t))]-E[\Psi(x)]=\int_{0}^{t}E[\mathcal{L}\Psi(X^{x}(s))]ds. 
    \end{align*}
     Then we substitute $[P_{t}\Psi](x)=E[\Psi(X^{x}(t))]$, and we get 
     \begin{align*}
     [P_{t}\Psi](x)-\Psi(x)=\int_{0}^{t}[P_{s}\mathcal{L}\Psi](x)ds.
     \end{align*}
     Again to prove (\ref{relation2}), we can rewrite (\ref{imldeq6t}) as,
     \begin{align*}
    d\Psi(X^{x}(t))=&\left\langle \frac{d\Psi(X^{x}(t))}{dx}, AX^{x}(t)+F(X^{x}(t))\right\rangle_{H}dt&   \\
    &+\frac{1}{2}tr\left(\frac{d^{2}\Psi(X^{x}(t))}{dx^{2}}(B(X^{x}(t))Q^{1/2})(B(X^{x}(t))Q^{1/2})^{*}\right)dt&  
    \end{align*}
    \begin{align*}
    +\int_{H\backslash\left\{0\right\}}\left[\Psi(X^{x}(t)+f(v,X^{x}(t)))-\Psi(X^{x}(t))-\left\langle \frac{d\Psi(X^{x}(t))}{dx},f(v,X^{x}(t))\right\rangle_{H}\right]\beta(dv)dt    
    \end{align*}
    \begin{align*}
    +\left\langle \frac{d\Psi(X^{x}(t))}{dx},B(X^{x}(t))dW_{t}\right\rangle_{H}	 
    \end{align*}
    \begin{align*}
    +\int_{H\backslash\left\{0\right\}}\left[\Psi(X^{x}(t)+f(v,X^{x}(t)))-\Psi(X^{x}(t))\right]q(dv,dt).	
    \end{align*}
    Now, since $[P_{t}\Psi](x)=E[\Psi(X^{x}(t))]$, so by Lebesgue dominated convergence theorem 
    \begin{align*}
    \lim_{t\downarrow0}\frac{[P_{t}\Psi](x)-\Psi(x)}{t}	
    \end{align*}
    \begin{align*}
    =E\lim_{t\downarrow0}\frac{1}{t}\int_{0}^{t}\left\langle \frac{d\Psi(X^{x}(s))}{dx}, AX^{x}(s)+F(X^{x}(s))\right\rangle_{H}ds	
    \end{align*}
    \begin{align*}
    +\frac{1}{2}E\lim_{t\downarrow0}\frac{1}{t}\int_{0}^{t}tr\left(\frac{d^{2}\Psi(X^{x}(s))}{dx^{2}}(B(X^{x}(s))Q^{1/2})(B(X^{x}(s))Q^{1/2})^{*}\right)ds	
    \end{align*}
    \begin{align*}
    +E\lim_{t\downarrow0}\frac{1}{t}\int_{0}^{t}\int_{H\backslash\left\{0\right\}}\left[\Psi(X^{x}(s)+f(v,X^{x}(s)))-\Psi(X^{x}(s))-\left\langle \frac{d\Psi(X^{x}(s))}{dx},f(v,X^{x}(s))\right\rangle_{H}\right]\beta(dv)ds 	
    \end{align*}
    \begin{align*}
    =\left\langle \frac{d\Psi(x)}{dx}, Ax+F(x)\right\rangle_{H}+\frac{1}{2}tr\left(\frac{d^{2}\Psi(x)}{dx^{2}}(B(x)Q^{1/2})(B(x)Q^{1/2})^{*}\right)	
    \end{align*}
    \begin{align*}
    +\int_{H\backslash\left\{0\right\}}\left[\Psi(x+f(v,x))-\Psi(x)-\left\langle \frac{d\Psi(x)}{dx},f(v,x)\right\rangle_{H}\right]\beta(dv)	
    \end{align*}
    \begin{align*}
    =[\mathcal{L}\Psi](x).
    \end{align*}
    Hence the proof follows.
  \end{proof}
  \end{theorem}


  Now we consider also the case when $X^x(t)\notin\mathcal{D}(A)$. $\Psi$ is independent of time $t$ and 
  $\Psi\in\tilde{C}_b^{2}(H)$, where $\tilde{C}_b^{2}(H)$ is defined like $\tilde{C}_b^{1,2}(H)$ however closing the operator $\mathcal{L}$ over functions defined on $C_b^2(H)$, i.e. $\tilde{C}_b^2(H)$ is  the class of functions 
  $\Psi\in C_b^{2}(H)$ satisfying  the properties: \\
 ($\tilde{1}$) The function $\mathcal{L}\Psi(x)$ can be extended to a continuous function $\overline{\mathcal{L}\Psi}(x)$ on $H$ for $x\in H$.\\
 ($\tilde{2}$) $\|\overline{\mathcal{L}\Psi}(x)\|
 \leq k(1+\|x\|^{2})$, for $x\in H$ and for some $k>0$.
   
   
  \begin{corollary}\label{ichigeneratorCor}
  Assume that $\overline{\mathcal{L}\Psi}$ satisfies the condition
  \begin{align}\label{ichiCont}
 \|\overline{\mathcal{L}\Psi}(y)-\overline{\mathcal{L}\Psi}(z)\|\leq k\|y-z\|^{2}(\|y\|^{2}+\|z\|^{2})	
 \end{align}
 for some $k>0$. Then 
 \begin{align}\label{exTrelation1}
  [P_{t}\Psi](x)-\Psi(x)=\int_{0}^{t}[P_{s}\overline{\mathcal{L}\Psi}](x)ds, \quad \text{for} \,\,\, x\in \mathcal{D}(A)
  \end{align}
  and 
  \begin{align}\label{exTrelation2}
  \lim_{t\downarrow0}\frac{[P_{t}\Psi](x)-\Psi(x)}{t}=[\overline{\mathcal{L}\Psi}](x), \quad \text{for} \,\,\, x\in \mathcal{D}(A).
  \end{align}
  \begin{proof} For the proof we refer to A. Ichikawa \cite{8}. It can be shortly described as follows:
   to prove \eqref{exTrelation1}, we take the expectation on both sides of \eqref{ichiIto1}  and use the relation 
   $[P_{t}\Psi](x)=E[\Psi(X^{x}(t))]$. The condition (\ref{ichiCont}) assures the continuity of 
   $P_s\overline{\mathcal{L}\Psi}$ in $s$. Therefore we can write
 \begin{align*}
 \lim_{t\downarrow0}\frac{[P_{t}\Psi](x)-\Psi(x)}{t}=E\lim_{t\downarrow0}\frac{1}{t}\int_0^t\overline{\mathcal{L}\Psi}(X^{x}(s))ds=[\overline{\mathcal{L}\Psi}](x).
 \end{align*}
 Hence the proof.
  \end{proof}
  \end{corollary}
  
  We define 
   
   \begin{align}\label{generatorDcor}
  [\mathcal{A}\Psi](x):=[\overline{\mathcal{L}\Psi}](x), \quad \text{for} \,\,\, x\in \mathcal{D}(A), \Psi\in\tilde{C}_b^{2}(H).
  \end{align}
  
 We call $\mathcal{A}$  the "weak generator of the Markov process $X^x(t)$".
 
 \vspace{0.5cm} 
  The existence of $\overline{\mathcal{L}\Psi}$ is rather restrictive. Therefore we introduce a class of functions which is larger than 
  $\tilde{C}_b^{2}(H)$, where continuous extension of the function $\mathcal{L}\Psi$ is not required. 


\begin{lemma}\label{newV}
Let $\Psi(x)$$\in C^{2}(H)$ satisfy the property: \\
(i) The function $\mathcal{L}\Psi(x)\leq\mathcal{U}(x)$, for $x\in\mathcal{D}(A)$, where $\mathcal{U}(x)$ is a continuous function on $H$ such that  $\|\mathcal{U}(x)\|\leq k(1+\|x\|^2)$ for some 
$k>0$.\\
 Moreover assume that conditions 
\eqref{bd-h-1}, \eqref{bd-h-2}, \eqref{bd-h-3} hold. 

Then
\begin{align}\label{newV1}
 \Psi(X^{x}(t))-\Psi(x)\leq&\int_0^t\mathcal{U}(X^{x}(s))ds+\int_{0}^{t}\left\langle \Psi_{x}(X^{x}(s)),B(X^{x}(s))dW_{s}\right\rangle_{H}&\nonumber\\
 &+\int_{0}^{t}\int_{H\backslash\left\{0\right\}}\left[\Psi(X^{x}(s)+f(v,X^{x}(s)))-\Psi(X^{x}(s))\right]q(dv,ds).&
 \end{align}
If, in particular, $\mathcal{U}(x)=0$, then $\Psi(X^{x}(t))$ is a supermartingale.
\begin{proof}
 $\int_0^t\mathcal{U}(X^{x}(s))ds$ is well defined. By Corollary \ref{itocorollary} 
\begin{align*}
\int_0^t\lim_{n\rightarrow\infty}\mathcal{L}_{n}\Psi(X_{n}^{x}(s))ds\leq\int_0^t\mathcal{U}(X^{x}(s))ds.
\end{align*}
Therefore we  conclude (\ref{newV1}).\\

For the second part, we take the conditional expectation in both sides of (\ref{newV1}), then we get

\begin{align*}
 E[\Psi(X^{x}(t))|\mathcal{F}_0^X]-E[\Psi(x)|\mathcal{F}_0^X]\leq0
 \end{align*}
 since $\mathcal{U}(x)=0$ and terms containing the Gaussian and non-Gaussian noise in the R.H.S. of \eqref{newV1} are martingales.
 \begin{align*}
 \Rightarrow E[\Psi(X^{x}(t))|\mathcal{F}_0^X]\leq\Psi(x).
\end{align*}
Therefore a supermartingale.
\end{proof}
\end{lemma}

\vspace{0.5cm}
 
 Now, through an example, we show how  Ichikawa type inequality \eqref{newV1} of Lemma \ref{newV} might be used to prove that 
 a mild solution is exponentially stable in the mean square sense.

\noindent\textbf{Example 4:-}\\
Consider the stochastic heat equation
\begin{align}\label{heat}
dX(x,t)=\frac{\partial^2}{\partial x^2}X(x,t)dt+B(X(x,t))dW_t+\int_{H\backslash\left\{0\right\}}f(v)X(x,t)q(dv,dt),
\end{align}
 with
\begin{align*}
0<x<1.
\end{align*}
\begin{align*}
 X(0,t)=X(1,t)=0;\ \ X(x,0)=X_0(x); \ \ X_0,B,f\in L_2(0,1).
\end{align*}

Here we take $H=L_2(0,1)$, $A=d^2/dx^2$ and 

\begin{align*}
\mathcal{D}(A)=\left\{f\in H|f', f''\in H; f(0)=f(1)=0\right\}.
\end{align*}

Since $A$ has eigenvectors $\left\{\sqrt{2}\sin n\pi x\right\}$ and eigenvalues $\left\{-n^2\pi^2\right\}$ for $n\in\mathbb{N}$. 
Then $X\in\mathcal{D}(A)$, $\left\langle AX,X\right\rangle\leq -\pi^2\|X\|^2$ (see example 6.1 of \cite{8}). \\

Now consider the function $\Psi(x)=\left\|x\right\|^2$. Therefore $\Psi_x(x)=2x$ and $\Psi_{xx}=2$. Hence, for $x\in\mathcal{D}(A)$

\begin{align*}
\left\langle \Psi_x(x),Ax\right\rangle\leq-2\pi^2\left\|x\right\|^2,
\end{align*}

\begin{align*}
\frac{1}{2}tr(\Psi_{xx}(x)(B(x))Q(B(x))^*)=tr((B(x))Q(B(x))^*)\leq l(1+\left\|x\right\|^2),
\end{align*}

by (A3) and 

\begin{align*}
&\int_{H\backslash\left\{0\right\}}[\Psi(x+f(v)x)-\Psi(x)-\left\langle \Psi_x(x),f(v)x\right\rangle]\beta(dv)&\\
&\leq\int_{H\backslash\left\{0\right\}}\left\|f(v)x\right\|^2\times\sup_{0\leq\theta\leq1}\|\Psi_{xx}(x+\theta f(v)y)\|\beta(dv)&\\
&\leq 2l(1+\left\|x\right\|^2)&\ \ \ \ [\text{since},\ \Psi_{xx}=2.]
\end{align*}

by (A3). To get the above inequality we followed the argument used in the proof of Theorem \ref{imldtm4}. Therefore for 
$x\in\mathcal{D}(A)$,

\begin{align*}
\mathcal{L}\Psi(x)=\mathcal{L}\left\|x\right\|^2&\leq-2\pi^2\left\|x\right\|^2+(l+2l)(1+\left\|x\right\|^2)&\\
&=(-2\pi^2+3l)\|x\|^2+3l&
\end{align*}
So that $\Psi$ satisfies (i)  in Lemma\ref{newV}.

\begin{align}\label{ineqVspc}
 \|X(t)\|^2-\|X(0)\|^2\leq&\int_0^t\{(-2\pi^2+3l)\|X(s)\|^2+3l\}ds&\\
 &+\int_0^t\langle2X(s),B(X(s))dW_s\rangle&\nonumber\\
 &+\int_0^t\int_{H\backslash\left\{0\right\}}[\|X(s)+f(v)X(s)\|^2-\|X(s)\|^2]q(dv,ds).\nonumber&
\end{align}

So, from Lemma \ref{newV}, we obtain the inequality of \eqref{ineqVspc} for the function $\Psi(x)=\left\|x\right\|^2$.\\

Now whenever $\pi^2>\frac{3}{2}l$, then for some constant $k>0$
\begin{align}\label{exp-V.m.s.s}
 \|X(t)\|^2-\|X(0)\|^2\leq&\int_0^t\{-k\|X(s)\|^2+3l\}ds
 +\int_0^t\langle2X(s),B(X(s))dW_s\rangle&\\
 &+\int_0^t\int_{H\backslash\left\{0\right\}}[\|X(s)+f(v)X(s)\|^2-\|X(s)\|^2]q(dv,ds).\nonumber&
\end{align}
 Now applying expectation on both sides of \eqref{exp-V.m.s.s}
 \begin{align}\label{exp-V-Nw1}
  E[\|X(t)\|^2]&\leq-k\int_0^tE[\|X(s)\|^2]ds+3lt+\|X(0)\|^2&\\
  &\leq-k\int_0^tE[\|X(s)\|^2]ds+\lambda\|X(0)\|^2\ \ \ \text{[for sufficiently large $\lambda>0$].}\nonumber&
 \end{align}
 Therefore, by Gronwall's lemma
 \begin{align*}
 E[\|X(t)\|^2]\leq\lambda e^{-kt}\|X(0)\|^2,
 \end{align*}
hence by Definition \ref{m.s.s-ex}, we  conclude that the solution $X(t)$ is exponentially stable in the mean square sense.

\section{ Mild It\^o formula written in terms of the semigroup}\label{pratojenroc} 
In \cite{11} an It\^o formula written in terms of the semigroup $\{S(t)\}_{t \geq 0}$ of the Drift operator $A$ was derived  for the Gaussian case. 
  There the mild It\^o process was transformed  to a standard It\^o process to which   the standard It\^o formula was applied. By relating this transformed standard It\^o process with
the original mild It\^o process with a suitable relation, the mild It\^o formula written in terms of  $\{S(t)\}_{t\geq 0}$ was obtained.\\
We present here the method and result for the non Gaussian case\\

 We consider the SPDE with values in $H$ as
 \begin{align}\label{withoutGauss}
  dX(t)=(AX(t)+F(X(t)))dt+\int_{H\backslash\left\{0\right\}} f(v,X(t))q(dv,dt);
 \end{align}  
 with initial condition $X(0)\in H$. Where $H$ is a real separable Hilbert space, $A$ is the generator of  $S(t)$  a pseudo-contraction semigroup $\{S(t)\}_{t\geq 0}$ and conditions (A1), (A3) are satisfied.
 (To simplify formulas we skip the  Gaussian term ). Then from 
 Definition \ref{imlddf1}, the mild solutions are defined as
 \begin{align*}
 X(t)=S(t)X(0)+\int_{0}^{t}S(t-s)F(X(s))ds+\int_{0}^{t}\int_{H\backslash\left\{0\right\}} S(t-s)f(v,X(s))q(dv,ds)
 \end{align*}
with probability one for $t\in[0,T]$. 
\\


\begin{theorem}\label{TheoremA}
Assume that $\Psi\in C^{1,2}([0,T]\times H)$, $\Psi:[0,T]\times H\rightarrow\mathbb R$.  Also 
assume that the conditions \eqref{bd-h-1}, \eqref{bd-h-2}, \eqref{bd-h-3} hold.  Then 
the following mild It\^o formula holds

\begin{align}\label{rocknerIto}
 &\Psi(t,X(t))=\Psi(S(t)X(0))+\int_{0}^{t}(\partial_{1}\Psi)(s,S(t-s)X(s))ds&\\
 &+\int_{0}^{t}(\partial_{2}\Psi)(s,S(t-s)X(s))S(t-s)F(X(s))ds&\nonumber\\ 
 &+\int_{0}^{t}\int_{H\backslash\left\{0\right\}}[\Psi(s,S(t-s)X(s)+S(t-s)f(v,X(s)))-\Psi(s,S(t-s)X(s))&\nonumber\\
 &\ \ \ \ \ \ \ \ \ \ \ \ \ \ \ \ \ \ -\left\langle(\partial_{2}\Psi)(s,S(t-s)X(s)),S(t-s)f(v,X(s))\right\rangle]\beta(dv)ds&\nonumber\\
 &+\int_{0}^{t}\int_{H\backslash\left\{0\right\}}\left[\Psi(s,S(t-s)X(s)+S(t-s)f(v,X(s)))-\Psi(s,S(t-s)X(s))\right]q(dv,ds)&\nonumber
 \end{align}
 
 $P$-a.s. for all $t\in[0,T]$. Here $(\partial_{1}\Psi)(t,x)=(\frac{\partial\Psi}{\partial t})(t,x)$ and $(\partial_{2}\Psi)(t,x)=(\frac{\partial\Psi}{\partial x})(t,x)$.              $(\partial_{1}\Psi)\in C([0,T]\times H,\mathbb{R})$ and $(\partial_{2}\Psi)\in C([0,T]\times H,\mathcal{L}(H,\mathbb{R}))$.

\begin{proof}
Here we use the transformation technique, given in \cite{14}. For existence and uniqueness of the mild solutions w.r.t. cPrm we refer \cite{4}.  \\
We remark that the integrals in \eqref{rocknerIto} are well defined due to the conditions \eqref{bd-h-1}, \eqref{bd-h-2}, \eqref{bd-h-3}. This can be proven in a similar way as in the It\^o formula for SDEs proven in Theorem 3.7.2 in \cite{bookR}.\\
 Let $U_t\in\mathcal{L}(H)$, $t\in[0,\infty)$, is a strongly continuous pseudo-contractive
  semigroup on $H$ and $S(t-s)=U_{(t-s)}\in\mathcal{L}(H)$ for all $0\leq s\leq t\leq T$.  \\
	
  Then there exists a separable $\mathbb{R}$-Hilbert space $(\mathcal{H},\langle.,.\rangle_{\mathcal{H}},\|.\|_{\mathcal{H}})$ with 
  $H\subset\mathcal{H}$ and $\|v\|_{H}=\|v\|_{\mathcal{H}}$ for all $v\in H$ and a strongly continuous group 
  $\mathcal{U}_t\in\mathcal{L}(\mathcal{H})$, $t\in\mathbb{R}$,(Here we use the fact that, the strongly continuous pseudo-contractive semigroup can be 
  dilated to strongly continuous group); such that 
  \begin{align}\label{trnsfrm}
  U_{t}(v)=P(\mathcal{U}_{t}(v))
 \end{align}
 for all $v\in H\subset\mathcal{H}$ and all $t\in[0,\infty)$ where $P:\mathcal{H}\rightarrow H$ is the orthogonal projection from 
 $\mathcal{H}$ to $H$.  \\

 Now we transform our mild It\^o process $X:[0,T]\times\Omega\rightarrow H$ to a standard It\^o process $\bar{X}$, by following the technique of 
 \cite{14}, roughly speaking by multiplying the mild It\^o process with $\mathcal{U}_{-t}$ for $t\in[0,T]$. 
Let, $\bar{X}:[0,T]\times\Omega\rightarrow\mathcal{H}$ be the unique adapted, c\`adl\`ag stochastic process such that,

\begin{align}\label{dar1}
 \bar{X}_{t}=X(0)+\int_{0}^{t}\mathcal{U}_{-s}F(X(s))ds+\int_{0}^{t}\int_{\mathcal{H}\backslash\left\{0\right\}}\mathcal{U}_{-s}f(v,X(s))q(dv,ds).
 \end{align}
 $P$-a.s. for all $t\in[0,T]$. Here we use the following transformation,

\begin{align}\label{dar2}
 P(\mathcal{U}_{t}(\bar{X}_{s}))&=P(\mathcal{U}_{t}(X(0)))+\int_{0}^{s}P\mathcal{U}_{(t-u)}F(X(u))du
 +\int_{0}^{s}\int_{H\backslash\left\{0\right\}}P\mathcal{U}_{(t-u)}f(v,X(u))q(dv,du)&  \\
 &=S(t)X(0)+\int_{0}^{s}S(t-u)F(X(u))du+\int_{0}^{s}\int_{H\backslash\left\{0\right\}}S(t-u)f(v,X(u))q(dv,du)&\nonumber  \\
 &=S(t-s)\left(S(s)X(0)+\int_{0}^{s}S(s-u)F(X(u))du+\int_{0}^{s}\int_{H\backslash\left\{0\right\}}S(s-u)f(v,X(u))q(dv,du)\right)&\nonumber  \\
 &=S(t-s)X(s).&\nonumber
 \end{align}
 \begin{align*}
 \Rightarrow P(\mathcal{U}_{t}(\bar{X}_{t}))=X(t)\ \ \text{and}\ \ P(\mathcal{U}_{t}(\bar{X}_{0}))=S(t)X(0).
 \end{align*}
 
 $P$-a.s. for all $s,t\in[0,T]$ with $s\leq t$.  \\

Now we will apply the It\^o formula for strong solutions of \cite{9} to the test function $\Psi(s,P(\mathcal{U}_{t}(v)))$ 
 for $s\in[0,t]$, $v\in\mathcal{H}$

\begin{align}\label{dar3}
 &\Psi(t,X(t))=\Psi(t,P(\mathcal{U}_{t}(\bar{X}_{t})))=\Psi(P(\mathcal{U}_{t}(\bar{X}_{0})))
       +\int_{0}^{t}(\partial_{1}\Psi)(s,P(\mathcal{U}_{t}(\bar{X}_{s})))ds&\\
 &+\int_{0}^{t}(\partial_{2}\Psi)(s,P(\mathcal{U}_{t}(\bar{X}_{s})))P\mathcal{U}_{(t-s)}F(X(s))ds&\nonumber\\ 
 &+\int_{0}^{t}\int_{H\backslash\left\{0\right\}}[\Psi(s,P(\mathcal{U}_{t}(\bar{X}_{s}))+P\mathcal{U}_{(t-s)}f(v,X(s)))-\Psi(s,P(\mathcal{U}_{t}(\bar{X}_{s})))&\nonumber\\
 &\ \ \ \ \ \ \ \ \ \ \ \ \ \ \ \ \ \ -\left\langle(\partial_{2}\Psi)(s,P(\mathcal{U}_{t}(\bar{X}_{s}))),P\mathcal{U}_{(t-s)}f(v,X(s))\right\rangle]\beta(dv)ds&\nonumber\\
 &+\int_{0}^{t}\int_{H\backslash\left\{0\right\}}\left[\Psi(s,P(\mathcal{U}_{t}(\bar{X}_{s}))+P\mathcal{U}_{(t-s)}f(v,X(s)))
      -\Psi(s,P(\mathcal{U}_{t}(\bar{X}_{s}))\right]q(dv,ds)&\nonumber
 \end{align}

$P$-a.s. for all $\Psi\in C^{1,2}([0,T]\times H)$. Now substituting \eqref{trnsfrm} and 
 (\ref{dar2}) in (\ref{dar3}), we get our Mild It\^o formula of (\ref{rocknerIto}).
 \end{proof}
\end{theorem}

  \begin{remark}
  The formula \eqref{rocex}  for the function $\Psi(s,x)=e^{(t-s)A}\varGamma(x)$ was obtained  through the It\^o formula \eqref{ito-lim-forml} in Theorem \ref{imldtm4}, as well as  the "Ichikawa type" -It\^o formula \eqref{ichiIto1} in Theorem \ref{ichiIto} and is obtained  also applying the It\^o formula \eqref{rocknerIto} in Theorem\ref{TheoremA}
  \end{remark}

\medskip
\noindent
{\bf Acknowledgements}: We are very grateful to Prof. P. Sundar (Louisiana State University) for useful discussions related to this work

\end{document}